\PassOptionsToPackage{colorlinks=true,urlcolor=blue,citecolor=red,
linkcolor=blue,linktocpage,pdfpagelabels,
bookmarksnumbered,bookmarksopen}{hyperref}

\documentclass[reqno,10pt]{amsart}
\usepackage{pgfplots}
\pgfplotsset{compat=1.18}
\usepackage{float}
\usepackage{enumitem}
\usepackage[dvipsnames,svgnames]{xcolor}
\usepackage[english]{babel}
\usepackage{csquotes}
\usepackage[left=2.8cm,right=2.8cm,top=2.9cm,bottom=2.9cm]{geometry}
\usepackage{graphicx,subcaption}
\usepackage{tikz}
\usetikzlibrary{arrows.meta}
\usepackage{calc}

\usepackage{amsthm}
\usepackage{mathtools}
\DeclarePairedDelimiter{\norm}{\lVert}{\rVert}

\newcommand{\de}{\mathrm{d}}
\usepackage{bbm}
\newtheorem{theorem}{Theorem}
\newtheorem{lemma}[theorem]{Lemma}

\newtheorem{proposition}[theorem]{Proposition}

\newtheorem{theoremletter}{Theorem}

\theoremstyle{definition}
\newtheorem{definition}[theorem]{Definition}
\theoremstyle{remark}
\theoremstyle{definition}

\newcommand{\innerthmname}{}

\theoremstyle{definition}

\makeatletter
\def\namedlabel#1#2{\begingroup
	#2
	\def\@currentlabel{#2}
	\phantomsection\label{#1}\endgroup
}
\makeatother

\DeclareMathOperator{\dist}{dist}

\DeclareMathOperator{\cat}{cat}

\DeclareMathOperator{\diam}{diam}

\usepackage[
  backend=biber,
  style=numeric,     % or authortitle
  sorting=nyt,          % sort by name, year, title
  giveninits=true,      % abbreviate first names
  doi=false,
  url=true,
  isbn=false,
  eprint=false,
  backref=true,
  maxnames=99
]{biblatex}
\DefineBibliographyStrings{english}{backrefpage = {Cited on page},backrefpages = {Cited on pages},}
\renewbibmacro{in:}{}
\DeclareFieldFormat[article]{title}{#1} 
\DeclareFieldFormat[article]{journaltitle}{\mkbibemph{#1}} 
\DeclareFieldFormat[article]{volume}{\textbf{#1}} 
\definecolor{TODOcolor}{RGB}{200,50,50}   % red
\definecolor{DONEcolor}{RGB}{40,150,40}    % green

\newif\ifdraft
\draftfalse   % disabled for arXiv submission

\ifdraft
  \usepackage[colorinlistoftodos,prependcaption,textsize=footnotesize]{todonotes}
\else
  \usepackage[disable]{todonotes} % same commands, but nothing is printed
\fi

\ifdraft
  \usepackage{pdfcomment} % clickable PDF annotations
\else
  \usepackage[final]{pdfcomment} % defines commands but hides notes
\fi

\title[Survey on topological methods for Allen--Cahn equations and systems]{Survey on topological methods for Allen--Cahn equations and systems}
\author[J.H. Andrade]{Jo\~{a}o Henrique Andrade*}
\author[S. Nardulli]{Stefano Nardulli}
\author[R. Ponciano]{Raon\'i Ponciano}

\address[J.H. Andrade]{Institute of Mathematics and Statistics,
	University of S\~ao Paulo
	\newline\indent 
05508-090, S\~ao Paulo-SP, Brazil}
\email{\href{mailto:andradejh@ime.usp.br}{andradejh@ime.usp.br}}

\address[S. Nardulli]{Department of Mathematics,
	Federal University of ABC
	\newline\indent 
09210-580, S\~ao Paulo-SP, Brazil}
\email{\href{mailto:stefano.nardulli@ufabc.edu.br}{stefano.nardulli@ufabc.edu.br}}

\address[R. Ponciano]{Department of Mathematics,
	Federal University of ABC
	\newline\indent 
09210-580, S\~ao Paulo-SP, Brazil}
\email{\href{mailto:raoni.ponciano@ufabc.edu.br}{raoni.ponciano@ufabc.edu.br}}

\thanks{* Corresponding author.}

\begin{document}

\begin{abstract}
We present a survey on multiplicity results for the Allen--Cahn equation and systems in the singular perturbation regime, emphasizing their geometric interpretation through $\Gamma$-convergence and isoperimetric theory. In the scalar case, the Allen--Cahn functional converges to perimeter, giving rise to minimal and constant-mean-curvature hypersurfaces, while vectorial Allen--Cahn systems lead to multi-phase isoperimetric clusters.
The main methodological tool discussed is the photography method, a variational-topological approach based on localized approximate solutions and barycenter maps, which enables one to encode the topology of the ambient manifold into multiplicity results. We compare problems posed on closed manifolds with those on manifolds with boundary, describing the distinct geometric effects induced by Neumann and Dirichlet boundary conditions.
The survey highlights both the effectiveness and the limitations of this framework, particularly in the vectorial case, where the lack of a full classification of isoperimetric clusters creates fundamental analytical challenges.
\end{abstract}

% \begin{center}
% {\textsl{Dedicated to Professor Paolo Piccione on the occasion of his 60th birthday}}
% \end{center}

\maketitle

\begin{center}
\footnotesize
\tableofcontents
\end{center}

%====================================================
%====================================================
\section{Introduction}
%====================================================
%====================================================
The Allen--Cahn equation plays a central role in the modern theory of phase transitions, serving as a diffuse-interface model that bridges nonlinear elliptic and parabolic equations with geometric evolution problems. In the singular limit, solutions develop sharp transition layers whose geometry is governed by classical isoperimetric problems and minimal hypersurfaces. This connection has made the Allen--Cahn equation a natural meeting point of nonlinear analysis, geometric measure theory, and global differential geometry.

A major theme of this survey is a variational-topological strategy, known as the \emph{photography method}, designed to produce multiple solutions of Allen--Cahn type problems in regimes where direct minimization or min-max constructions are difficult to implement. Roughly speaking, this method associates to each point of the underlying manifold a sharply localized approximate solution and uses barycenter maps together with Lusternik--Schnirelmann or Morse theory to transfer the topology of the ambient space into multiplicity results for the equation. This approach has proven particularly effective for problems with constraints, boundary effects, and vectorial structures.
    
The purpose of this survey is not to provide an exhaustive account of all developments related to the Allen--Cahn equation, but rather to present a coherent survey of this variational-topological framework, emphasizing both its scope and its limitations. Special attention is given to the contrast between problems posed on closed manifolds and those with boundary, and to the additional difficulties arising in vectorial systems.

In this introduction, we first provide a historical overview of the Allen--Cahn equation and its physical formulation. We then emphasize its significance by highlighting the many connections it establishes between various mathematical problems. Next, we present the Allen--Cahn equation in its full generality, considering its extension to systems, Riemannian manifolds, and both Dirichlet and Neumann boundary conditions. Finally, we summarize the main multiplicity results obtained via the photography method.

%====================================================
\subsection{Historical background}
%====================================================
The Allen--Cahn equation has its origins in early attempts to model phase transitions in materials with non-uniform composition. It traces back to the pioneering work of Van der Waals \cite{vanderwaals} in 1893 (see also the translated article \cite{MR523642}), where a diffuse interface approach was introduced to describe the coexistence of phases. In 1958, his idea resurfaced in the context of phase separation through the work of Cahn and Hilliard \cite{CAHN1961795,10.1063/1.1744102,10.1063/1.1730447}, who provided a mathematical formulation for the dynamics of such systems, now known as the Cahn--Hilliard equation. 

In 1978, Allen and Cahn \cite{ALLEN19791085} proposed a simplified evolution equation to describe the motion of interfaces in crystalline solids, highlighting an intrinsic connection between interface motion and mean curvature. This geometric insight laid the groundwork for a rich interplay between partial differential equations and geometric analysis, making the Allen--Cahn equation a central object in the study of phase transition phenomena.
From the modeling viewpoint, \emph{systems} arise naturally when more than two stable phases are
present, leading to multi-phase partitioning in the sharp-interface limit.

Over the years, the Allen--Cahn equation and its applications have been widely used across several scientific domains such as image inpainting \cite{LI201565}, two-phase incompressible fluids \cite{Ma2017867}, complex dynamics of dendritic growth \cite{TAKAKI201321,https://doi.org/10.1002/nme.5372}, mean curvature flows \cite{10.1093/imanum/dru058,Feng200333,Lee2015535}, image segmentation \cite{BENES2004187,LI201432}, and animal spots \cite{MR1908418,MR1952568}. Similarly, the Cahn–Hilliard equation has played a central role in modeling a wide variety of physical and biological phenomena, including model segregation-driven processes \cite{https://doi.org/10.1155/2016/9532608}, diblock copolymers \cite{JEONG20141263,doi:10.1021/ma00164a028}, topology optimization \cite{Zhou200789}, image inpainting \cite{doi:10.1137/060660631,4032803}, tumor growth \cite{doi:10.1142/S0218202516500263,doi:10.1142/S0218202510004313,WISE2008524}, biological phenomena \cite{PhysRevE.77.051129}, two-phase flows \cite{doi:10.1142/S0218202511500138,annurev:/content/journals/10.1146/annurev.fluid.30.1.139,doi:10.1142/S0218202596000341,RevModPhys.49.435,LAM_WU_2018}, moving contact line dynamics \cite{JACQMIN_2000,QIAN_WANG_SHENG_2006}, and multiphase incompressible fluid flow \cite{BADALASSI2003371,kotschote2015strong,LI201684,LI20161}.

%====================================================
\subsection{Phase separation and transitions}
%====================================================
The physical intuition behind the Allen--Cahn equation originates
from the modeling of phase separation processes in binary alloys,
that is, materials composed of two distinct components whose concentrations may vary in space.
Such systems tend to segregate into regions occupied by
energetically preferred phases, separated by thin transition layers.

\begin{figure}[ht]
    \centering
    \begin{tikzpicture}
	\begin{axis}[
		width=12cm,
		height=6cm,
		axis x line=middle,
		axis y line=left,
		xlabel={$x$},
		xlabel style={
			at={(ticklabel cs:1)},
			anchor=west,
		},
		ylabel={},
		xmin=2, xmax=10,
		ymin=-1.4, ymax=1.5,
		xtick=\empty,
		ytick={-1,1},
		tick label style={font=\small},
		label style={font=\small},
		x axis line style={-Latex,line width=0.6pt},
		y axis line style={-Latex,line width=0.6pt},
		]

% punti di "interfaccia"
\def\xa{4.5}
\def\xb{7.5}

% lunghezza tratto tratteggiato agli estremi
\def\delta{2}

% --- SINISTRA: plateau u = -1 (tratteggiato)
\addplot[
  blue!75!black,
  dashed,
  line width=0.9pt,
]
coordinates {(0,-1) (\xa -\delta,-1)};

% --- CENTRO: transizione liscia C^2
\addplot[
  blue!75!black,
  line width=0.9pt,
  domain=\xa-\delta:\xb + \delta,
  samples=400,
]
{ (x<=\xa) * (-1) + (x>=\xb) * ( 1) + ((x>\xa) * (x<\xb)) * ( -1 + 2 * ( 6*pow((x-\xa)/(\xb-\xa),5) -15*pow((x-\xa)/(\xb-\xa),4) +10*pow((x-\xa)/(\xb-\xa),3) ) ) };

% --- DESTRA: plateau u = 1 (tratteggiato)
\addplot[
  blue!75!black,
  dashed,
  line width=0.9pt,
]
coordinates {(\xb + \delta,1) (12,1)};

		% linee tratteggiate
		\addplot[dashed,line width=0.5pt]
			coordinates {(\xa,-1.4) (\xa,1.4)};
		\addplot[dashed,line width=0.5pt]
			coordinates {(\xb,-1.4) (\xb,1.4)};

		% etichette
		\node[black] at (axis cs:3.2,-1.3) { Phase A};
		\node[black] at (axis cs:6,-1.3) { Interface};
		\node[black] at (axis cs:8.8,-1.3) { Phase B};
		\node[black] at (axis cs:2.6,1.4) {$u(t,x)$};
	\end{axis}
    \end{tikzpicture}
    \caption{Phase transition, where $\{u=-1\}$ indicates the material is entirely in Phase A, while $\{u=1\}$ corresponds to the material being fully in Phase B.}
    \label{figure1}
\end{figure}
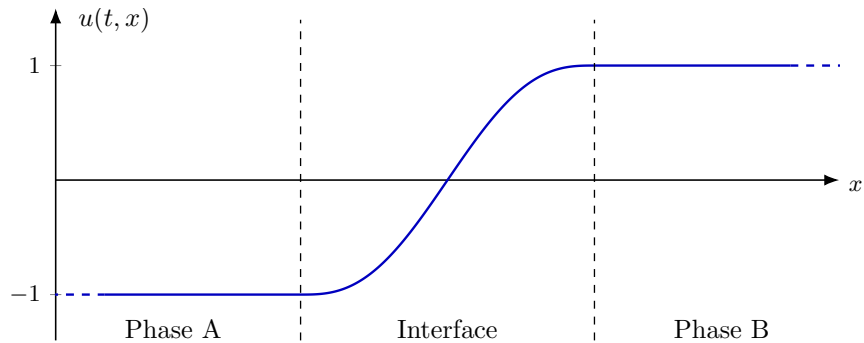

A convenient way to describe this behavior is through a continuous
scalar field, which encodes the local state of the system.
More precisely, let $\Omega \subset \mathbb{R}^N$ be a region, and let $u:\Omega\times \mathbb{R}_{\ge 0}\to \mathbb{R}$ represent the local state of the system at point $x\in \Omega$ and time $t\in\mathbb{R}_{\ge 0}$.
In the simplest and most classical setting of two competing phases,
this description is associated with a double-well potential
$W\colon\mathbb{R}\to\mathbb{R}$, whose minima represent the two pure phases.
Accordingly, the values $\{u=-1\}$ and $\{u=1\}$ correspond,
by a convenient normalization, to regions entirely occupied by Phase~A and Phase~B, respectively.
Phase transitions occur across interfaces, where the field $u$
varies smoothly between these two values
(see Figure~\ref{figure1}).

In this framework, the evolution of the system is governed by the
time-dependent (parabolic) Allen--Cahn equation, introduced in the
seminal work of Allen and Cahn~\cite{ALLEN19791085},
\begin{equation}\label{t-AC-PDE}
\tag{$t$-$AC_{N,\varepsilon,\Omega}$}
\dfrac{\partial u}{\partial t}
=
\varepsilon \Delta u - \frac{1}{\varepsilon}W'(u),
\quad \text{in } \Omega \times (0,\infty),
\end{equation}
where $\Delta$ denotes the Laplace operator with respect to the spatial variable $x\in\Omega$, $\varepsilon>0$ is a small parameter
related to the interfacial energy, and $W$ is the double-well potential to be discussed later.
In particular, $\varepsilon$ controls the thickness of the transition
layer between phases: smaller values of $\varepsilon$ correspond to sharper interfaces.
The stationary version of \eqref{t-AC-PDE}, which is also often referred to simply as the Allen--Cahn equation, is
\begin{equation}\label{AllenCahn}
\tag{$AC_{N,\varepsilon,\Omega}$}
\varepsilon\Delta u-\dfrac{1}{\varepsilon}W^{\prime}(u)=0,\quad\mbox{in }\Omega.
\end{equation}
Here, we have introduced a normalization (specifically, division by $\varepsilon$) as this scaling is essential later for the $\Gamma$-convergence.

\begin{figure}[ht]
    \centering
    \includegraphics[width=0.7\textwidth]{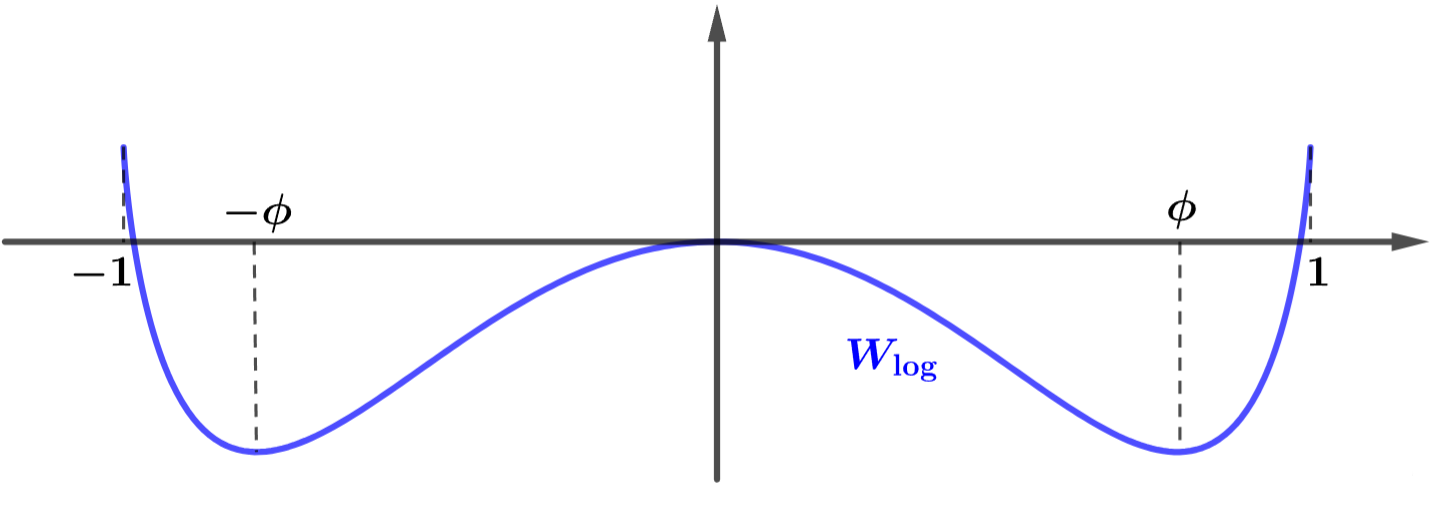}
    \caption{Graphic representation of $W_{\mathrm{log}}$ with $T=3$ and $T_c=4$.}
    \label{figure2}
    \end{figure}

While the specific choice of the values $\pm1$ is a convenient normalization,
the presence of a double-well structure arises naturally from thermodynamically consistent models.
In general terms, a double-well potential is any $C^2$ function $W$ that has exactly two local minima and is, in most cases, coercive (that is, $W(s)\to\infty$ as $s\to\pm\infty$). The first double-well potential was introduced by Cahn and Hilliard, motivated by thermodynamic principles (see \cite[Equation (3.1)]{10.1063/1.1744102} for the classical work and \cite[Section 4.2]{NOVICKCOHEN2008201} for a modern treatment). This original potential, often referred to as the logarithmic double-well potential, is given by
\begin{equation}\label{Wlog}
W_{\mathrm{log}}(s)=\dfrac{T}{2}\left[(1+s)\log(1+s)+(1-s)\log(1-s)\right]-\dfrac{T_c}{2}s^2,\quad s\in(-1,1)
\end{equation}
where $T>0$ denotes the absolute temperature of the system, and $T_c$ is the critical temperature at which phase separation occurs. When $0<T<T_c$, it is straightforward to verify that $W_{\mathrm{log}}$ exhibits a double-well structure, with two minima located at $\pm\phi\in(-1, 1)$, where $\phi$ is the positive root of the equation $W^{\prime}_{\mathrm{log}}(s)=0$ (see Figure \ref{figure2} for a plot of $W_{\mathrm{log}}$).

\begin{figure}[ht]
    \centering
    \includegraphics[width=0.5\textwidth]{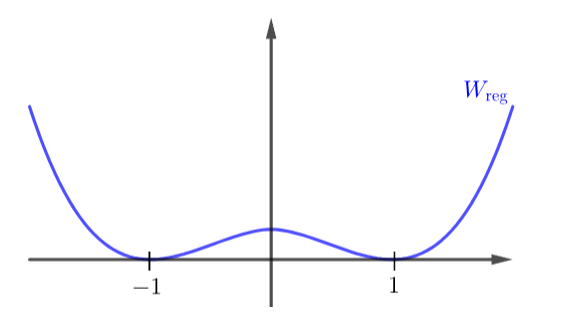}
    \caption{Graphic representation of $W_{\mathrm{reg}} = \frac{1}{4}(1-s^2)^2$}
    \label{figure3}
    \end{figure}

However, the presence of logarithmic terms in $W_{\mathrm{log}}$ makes analytical treatment difficult (cf.~\cite{MR4151195}). A common approach is to approximate $W_{\mathrm{log}}$ using a Taylor expansion. Around $s=0$, we have
\begin{equation*}
W_{\mathrm{log}}(s)=\dfrac{T-T_c}2s^2+\dfrac{T}{12}s^4+O(s^6) \quad {\rm as} \quad s\to 0.
\end{equation*}
Assuming $T=3$ and $T_c=4$, and neglecting higher-order terms and adding a constant to simplify the expression (without affecting the dynamics), we obtain the widely used regularized double-well potential:
\begin{equation}\label{Wreg}
W_{\mathrm{reg}}(s)=\frac14(1-s^2)^2.
\end{equation}
This regular potential $W_{\mathrm{reg}}$,
represented in Figure~\ref{figure3}, is employed in the majority of studies involving the Allen--Cahn equation due to its simplicity. 
Note, however, that this approximation modifies the location of the minima: while the original logarithmic potential has minima at $\pm\phi\in(-1,1)$, the regularized potential places them exactly at $\pm1$. From a physical perspective, this corresponds to assuming that the pure phases are stable, neglecting the presence of impurities in real materials.
For further insight into the physical interpretation of the Allen--Cahn equation, we refer the reader to the modern treatments presented in \cite{MR3309171,emmerich2003diffuse}, as well as to the classical work by \cite{ALLEN19791085}.

Another fundamental equation used to model phase transitions is the Cahn--Hilliard equation:
\begin{equation}\label{CahnHilliard}
\tag{$t$-$CH_{N,\varepsilon,\Omega}$}
\dfrac{\partial u}{\partial t}=\mathrm{div}\left[\mu(u)\nabla({\varepsilon}^{-1}W'(u)-\varepsilon\Delta u)\right],\quad\mbox{in }\Omega\times(0,\infty),
\end{equation}
where $\mathrm{div}$ and $\nabla$ denote the divergence and gradient with respect to the spatial variable $x\in\Omega$. The function $\mu\in C^2(\mathbb{R},\mathbb{R})$ represents the diffusion mobility, which may be taken either as a positive constant or as a concentration-dependent function (for more details, we refer \cite[Section 4]{NOVICKCOHEN2008201}). 
A comprehensive treatment of the Cahn--Hilliard equation can be found in \cite{MR4001523}. Although the Allen--Cahn equation and Cahn--Hilliard equation share structural similarities, neither is a strict simplification of the other.
A key distinction lies in their order and variational frameworks: while~\eqref{CahnHilliard} is a fourth-order PDE with a gradient flow in $H^{-1}$ topology and preserves mass,~\eqref{AllenCahn} is a second-order PDE with a gradient flow in $L^2$ and does not conserve mass.
For a detailed comparison between these two equations, we refer the reader to \cite{LI2019311}.

While the scalar Allen--Cahn equation
just discussed captures the emergence of interfaces
between two pure phases, it is natural to ask whether the same framework
can be extended to describe systems with more than two coexisting phases.
Classical examples include multiphasic fluid mixtures, polycrystalline
materials, and multi-component alloys, where several distinct pure states may appear and
interact through interfaces with different surface tensions.
One might attempt to introduce a scalar potential
$W\colon\mathbb R\to\mathbb R_{\geq0}$ with three or more wells. However, this approach is intrinsically
limited. Indeed, if $W$ has three distinct minima $z_0<z_1<z_2$, then any transition
from $z_0$ to $z_2$ must necessarily pass through the intermediate well $z_1$.
As a consequence, the resulting interface decomposes into two successive two-phase
transitions, and no genuinely new interface geometry is created. In particular, scalar
multi-well potentials do not produce true multiphase junctions in the singular limit,
and all interfaces reduce to combinations of two-phase boundaries.

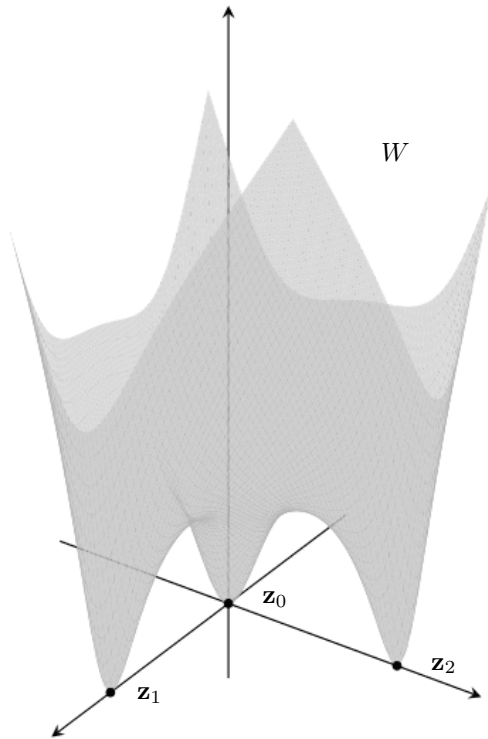
\begin{figure}
		\centering
		\begin{tikzpicture}
			\node[anchor=south west,inner sep=0] (image) at (0,0) {\includegraphics[width=0.4\textwidth]{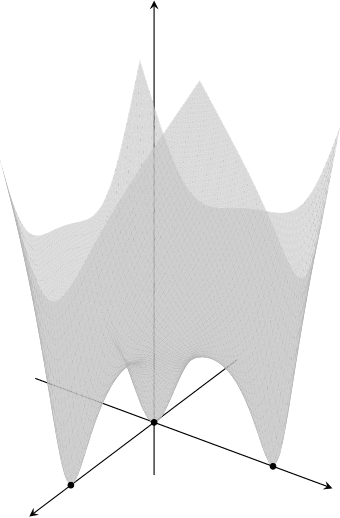}};
			\begin{scope}[x={(image.south east)},y={(image.north west)}]
				\node at (0.55,0.19) {$\mathbf{z}_0$};
				\node at (0.29,0.06) {$\mathbf{z}_1$};
				\node at (0.9,0.1) {$\mathbf{z}_2$};
				\node at (0.8,0.8) {$W$};
			\end{scope}
		\end{tikzpicture}
		\caption{
			A schematic representation of a multi-well potential in the vectorial setting,
illustrated here in the case of three phases, with noncollinear minima corresponding
to different pure states.
}
		\label{fig:potential-well}
	\end{figure}

To model genuine multiphase separation, it is therefore necessary to move beyond scalar
order parameters and consider vector-valued maps
$u\in C^2(\Omega,\mathbb R^m)$, together with a multi-well potential
$W\colon\mathbb R^m\to\mathbb R_{\geq0}$ whose minima set ${\rm argmin}(W)\subset \mathbb{R}^m$ consists of distinct and noncollinear points.
In this vectorial setting, the wells represent different phases that are
no longer totally ordered, and the geometry of transitions between phases becomes
essentially multidimensional.
A schematic example of such a multi-well structure is shown in
Figure~\ref{fig:potential-well}.
In particular, one can model up to $m+1$ distinct phases
using an $\mathbb R^m$-valued order parameter, and vectorial Allen--Cahn systems provide
a natural diffuse-interface approximation of multiphase problems.
At the level of diffuse-interface models, this leads to vectorial
parabolic Allen--Cahn systems, where the unknown is a vector-valued map
$u\colon\Omega\times [0,+\infty) \to\mathbb R^m$ satisfying
\begin{equation}\label{vAC}
\tag{$t$-$AC_{N,m,\varepsilon,\Omega}$}
\dfrac{\partial u}{\partial t}
=
\varepsilon \Delta u - \frac{1}{\varepsilon}\nabla W(u),
\quad \text{in } \Omega \times (0,\infty),
\end{equation}
where $W\in C^2(\mathbb R^m,\mathbb R_{\geq0})$ is a multi-well potential.

%====================================================
\subsection{Variational framework}
%====================================================
The Allen--Cahn equation admits a natural variational structure.
More precisely, the energy functional associated with~\eqref{AllenCahn} is given by
\begin{equation*}
\mathcal{AC}_\varepsilon(u)=
\left\{
\begin{array}{ll}
\displaystyle
\int_\Omega\left(
\dfrac{\varepsilon}{2}|\nabla u|^2+\dfrac{1}{\varepsilon} W(u)
\right)\,\mathrm dx,
& \mbox{if } u\in W^{1,2}(\Omega),\\[0.3cm]
\infty,
& \mbox{if } u\in L^1(\Omega)\setminus W^{1,2}(\Omega),
\end{array}
\right.
\end{equation*}
where the gradient term penalizes spatial variations of the field,
while the potential term penalizes deviations from the pure phases.

From a modeling perspective, the Allen--Cahn energy captures the leading-order
interfacial effects.
More refined descriptions may incorporate nonlocal interaction terms of Kac-type
to account for long-range interactions, leading to nonlocal phase-field models
with applications to the Ising process and to materials with microstructure
(see, for instance, \cite{MR3748585,MR1453735,MR1638739}).
In this survey, however, we focus on the local Allen--Cahn framework,
which already exhibits a rich interplay between variational structure,
geometry, and topology.

For small $\varepsilon>0$, the potential term $W(u)/\varepsilon$ strongly penalizes intermediate values of $u$, forcing low-energy configurations to stay
close to the pure phases $u\approx\pm1$ almost everywhere. As a consequence, any family of functions with uniformly bounded energy develops large regions where
$u^\varepsilon\approx+1$ and $u^\varepsilon\approx-1$, separated by a thin transition
layer whose thickness is proportional to $\varepsilon$. Inside these narrow
regions, the gradient term and the potential term balance each other, producing a steep but smooth transition between the
two phases.

From a geometric viewpoint, these transition layers behave as diffuse
hypersurfaces. As $\varepsilon\to0$, their thickness collapses to zero while their
location converges to a sharp interface $\Sigma\subset \Omega$.
A fundamental heuristic principle
is that the energy stored in such a layer is proportional, up to a universal
constant, to the $(N-1)$-dimensional measure of $\Sigma$. In other words, in the
singular limit, the Allen--Cahn energy counts how much interface is present, and
therefore approximates a perimeter functional. This is the reason why phase-field
models such as Allen--Cahn naturally select interfaces of minimal area, and why
their dynamics are closely related to mean curvature flow, while stationary
configurations correspond to constant mean curvature hypersurfaces under suitable
constraints.

This heuristic picture is made rigorous by the theory of
$\Gamma$-convergence, which shows that the Allen--Cahn energy $\mathcal{AC}_\varepsilon$
converges, as $\varepsilon\to0$, to the perimeter
functional of the limiting interface. A precise variational formulation of this
convergence, together with its extensions to manifolds, boundary conditions, and
vectorial systems, is presented in
Section~\ref{subsec:gamma_convergence}.

%====================================================
\subsection{Related geometric problems}
%====================================================
To emphasize the importance of the Allen--Cahn equation, we discuss its connection to different types of geometric-variational objects, such as the mean curvature flow (MCF), the constant mean curvature (CMC) hypersurfaces, and minimal partitions. 
\smallskip

\begin{figure}[ht]
    \centering
    \includegraphics[width=0.5\textwidth]{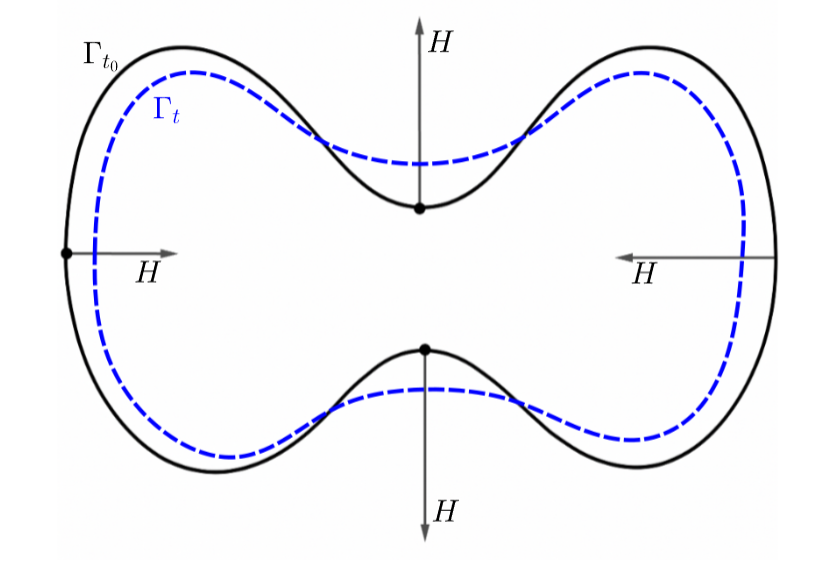}
    \caption{An interface following the mean curvature flow.}
    \label{figure4}
    \end{figure}

\noindent{\it Mean curvature flow.}
Let us consider the Allen--Cahn evolution equation with the regular double-well potential:
\begin{equation}
\tag{$t$-$AC_{N,\varepsilon}$}
\dfrac{\partial u}{\partial t}=\varepsilon^2\Delta u-W^{\prime}(u),\quad\mbox{in }\mathbb R^N\times(0,\infty),
\end{equation}
where the potential $W$ is given by \eqref{Wreg}. The gradient flow of the energy functional associated with this equation is
\begin{equation*}
\mathcal{AC}_\varepsilon(u(\cdot,t))=\int_{\mathbb R^N}\left(\dfrac{\varepsilon}2|\nabla u(x,t)|^2+\dfrac1\varepsilon W(u(x,t))\right)\mathrm dx,\quad u(\cdot,t)\in W^{1,2}(\mathbb R^N).
\end{equation*}
This transition layer defines an interface $\Gamma_t$, which evolves in time $t$ according to mean curvature flow (see Figure \ref{figure4}). We refer the reader to \cite{MR1205984,MR1101239,MR1153311,MR1100211,MR1177477,MR1237490,MR1341031,MR1055457,MR978829} for a rigorous formulation and proofs of these statements, and we cite \cite{MR1674799} for a more comprehensive list of references. 

For any fixed $t_0\in(0,\infty)$, the evolution of $\Gamma_{t_0}$ by the mean curvature flow is a smooth one-parameter family of immersions $\varphi\colon\Gamma_{t_0}\times[t_0,\infty)\to\mathbb R^{N}$ satisfying the system
\begin{equation*}
\left\{\begin{array}{ll}
     \dfrac{\partial}{\partial t}\varphi(x,t)=H(\varphi(x,t)),&\mbox{in }\Gamma_{t_0}\times[t_0,\infty)  \\
     \varphi(x,t_0)=x, 
\end{array}\right.
\end{equation*}
where $H(\varphi(x,t))$ denotes the mean curvature vector of the hypersurface $\Gamma_t=\varphi(\Gamma_{t_0},t)$ on the point $\varphi(x,t)$. We also refer to \cite{MR2815949} for lecture notes on mean curvature flow. 

The physical interpretation of the interfaces $\Gamma_t$ evolving through the mean curvature flow is particularly interesting: regions of Phase A that are surrounded by Phase B tend to transform into Phase B (and vice versa) due to surface tension effects. This transformation occurs more rapidly at points where the interface has high curvature, leading to an interface motion governed by mean curvature flow. The mean curvature vector indicates the direction in which the interface bends, and its magnitude determines the speed of the interface motion.

For the purpose of this survey, we focus on the time-independent case of the Allen--Cahn equation \eqref{AllenCahn} since it plays a fundamental role in the study of minimal surfaces. Indeed, stationary solutions of Eq. \eqref{t-AC-PDE} give rise to interfaces with zero mean curvature (that is, minimal surfaces). Moreover, due to the properties of the recovery sequences in the context of $\Gamma$-convergence, given a minimal surface, there exists a sequence of solutions $u^\varepsilon$ to Eq. \eqref{AllenCahn} such that the interface coincides with the minimal surface in the limit as $\varepsilon\to0$. This idea is central to the photography method, which we will describe in Section \ref{sec:photography}. The connection between the Allen--Cahn equation and minimal surfaces has been employed in the study of the famous De Giorgi's Conjecture; we refer the reader to the survey \cite{MR3864761} for a comprehensive overview. 
\smallskip

\noindent{\it From CMC hypersurfaces to minimal partitions.}
The relationship between the Allen--Cahn equation and minimal hypersurfaces places phase-field methods at the crossroads of geometric analysis and variational geometry. From this
perspective, the classical theory of minimal hypersurfaces naturally appears as the
zero-volume-constraint limit of diffuse interface models.

A central motivation for extending this framework comes from longstanding problems in
geometry concerning the existence and multiplicity of critical points of the area
functional under additional constraints. A prominent example is Yau's conjecture
\cite{MR645762}, which predicts the abundance and density of minimal hypersurfaces on closed
Riemannian manifolds. This conjecture has recently been resolved in several settings
\cite{MR3674223,MR3953507,MR4564260}, at least generically or in low dimensions.

Beyond minimal hypersurfaces, one is led to consider \emph{isoperimetric sets} and, more
generally, \emph{constant mean curvature} (CMC) hypersurfaces, which arise as critical
points of the perimeter functional under fixed-volume constraints. The situation becomes
substantially more delicate in this regime: while existence results are available for
isoperimetric regions \cite{MR4649390,MR4011704}, considerably fewer results are known for
\emph{multi-isoperimetric clusters}, where several volumes are prescribed simultaneously
\cite{MR4459029,MR2976521}. In contrast, the corresponding nonlocal isoperimetric problems
are by now much better understood; see, for instance,
\cite{arXiv:2306.07100} and the references therein.

From the diffuse-interface viewpoint, these geometric problems are naturally approximated
by Allen--Cahn type energies with suitable constraints. Indeed, on closed manifolds, the
Allen--Cahn energy $\Gamma$-converges, as $\varepsilon\to0$, to the perimeter functional,
and volume constraints at the PDE level give rise to constant mean curvature interfaces in
the sharp-interface limit. This observation underlies a PDE-based approach to geometric
variational problems developed in
\cite{MR4498838,MR3945835,MR3743704}, among others, and provides an alternative to the
classical tools of geometric measure theory.

These ideas extend naturally to multiphase settings. In this case, vectorial Allen--Cahn
systems approximate partitions of the ambient space into several regions, with interfaces
meeting along junctions governed by surface tension balance laws. The lack of a complete
classification of minimizers for the multi-isoperimetric problem makes the analysis
considerably harder, and it is precisely in this context that variational-topological
constructions, such as the photography method described later in this survey, become
particularly effective.

%====================================================
\subsection{Analytical formulation}
%====================================================
Throughout the paper, we denote by $m\in\mathbb Z_{>0}$ the number of chambers
(hence, for a two-phase system, $m = 1$). Also, let $(M,\partial M,g)$ denote an $N$-dimensional Riemannian manifold,
with or without boundary. Additionally, let $X\in\mathfrak X(M)$ be a smooth tangential vector field. We denote by $\nu$ the unit \emph{inner} normal vector field along $\partial M$. For a map $u=(u^1,\dots,u^m)\colon M\to\mathbb{R}^m$,
we define the vectorial gradient and Laplace--Beltrami operator componentwise,
\[
\nabla_g u := (\nabla_g u^1,\dots,\nabla_g u^m),
\qquad
\Delta_g u := (\Delta_g u^1,\dots,\Delta_g u^m).
\]
We work in the Sobolev space
\[
W^{1,2}(M,\mathbb{R}^m)
:=
\{u=(u^1,\dots,u^m): u^i\in W^{1,2}(M,\mathbb{R}),\ i=1,\dots,m\},
\]
equipped with the norm
\[
\|u\|_{W^{1,2}(M,\mathbb{R}^m)}
=
\Big(\sum_{i=1}^m \|u^i\|^2_{W^{1,2}(M,\mathbb{R})}\Big)^{1/2}.
\]

Now, let us precisely state the definition of the multi-well potential that will be used in the rest of the paper. 
\begin{definition}\label{def:multi-well}
 A nonnegative function $W\in C^2(\mathbb{R}^m,\mathbb{R}_{\geq0})$ is called a \emph{multi-well potential} if it satisfies the following conditions:
		\begin{enumerate}
			\item[\namedlabel{item:W-global-minima}{($W_0$)}] It has exactly $m + 1$ vanishing and nondegenerate global minima in $\mathbb{R}^m$, denoted by $\mathcal{Z}\coloneqq\{\mathbf{z}_{0},\dots,\mathbf{z}_{m}\}\subset\mathbb{R}^m$;
			this means that, for each $i=0,\ldots,m$, we have $W(\mathbf{z}_i)=0$, $\nabla W(\mathbf{z}_i)=0$, and $\nabla^2W(\mathbf{z}_i)>0$.
			We assume that $\mathbf{z}_0=0$,
			and $\{\mathbf{z}_{1},\dots,\mathbf{z}_{m}\}$
			is linearly independent.
			Moreover, the following \emph{strict} inequality holds:
			\begin{equation*}
				\tag{${{\rm W}_0}$}\label{eq:W0}
				\omega_{ij} <
				\omega_{i\ell} + \omega_{\ell j},
				\quad \forall i,j,\ell\in\{0,\dots, m\},\; \ell\notin\{i,j\},
			\end{equation*}
			where the $\omega_{ij}$-s are defined as follows:
			\begin{equation}
				\label{eq:def-omega-ij}
				\omega_{ij} 
				\coloneqq
				\inf_{\gamma\in \Gamma_{ij}}
				\int_{0}^{1}W^{1/2}(\gamma(s))\norm{\gamma^{\prime}(s)}\de s,
			\end{equation}
			and
			\begin{equation*}
				\Gamma_{ij} = \left\{
				\gamma \in C^1([0,1],\mathbb{R}^m):
				\gamma(0)=\mathbf{z}_i, \; \gamma(1)=\mathbf{z}_j
				\right\};
			\end{equation*} 
			\item[\namedlabel{item:W-def-k1}{{($W_1$)}}]
			There exists a constant $k_1>0$ such that 
			\begin{equation*}
				\tag{${{\rm W}_1}$}\label{eq:W1}
				|\nabla W(z)| \le k_1(1+|z|^{p-1}),
				\quad \forall z\in\mathbb{R}^m,
			\end{equation*}
			for some $1<p<2^{*}$ if $N\ge 3$ (or $p<\infty$ if $N=2$), where $2^{*}\coloneqq\frac{2N}{N-2}$ is the critical Sobolev exponent of the embedding $W^{1,2}(\mathbb R^N)\hookrightarrow L^q(\mathbb R^N)$ for $q>1$;
			\item[\namedlabel{item:W-def-k2}{{($W_2$)}}]
			There exists $k_2>0$ such that 
			\begin{equation}
				\tag{${{\rm W}_2}$}\label{eq:W2}
				|\nabla^2 W(z)|\le k_2(1+|z|^{p-2}),
				\quad \forall z\in\mathbb{R}^m
			\end{equation}
			for some $1<p<2^{*}$ if $N\ge 3$ (or $1<p<\infty$ if $N=2$);
			\item[\namedlabel{item:W-def-p12-k34-R}{($W_3$)}]
			There exist $p_1, p_2, k_3, k_4,R>0$ such that
			\begin{equation}\tag{${{\rm W}_3}$}
				\label{eq:W3}
				k_3|z|^{p_1}<W(z)<k_4|z|^{p_2},
				\quad \forall z\in\mathbb{R}^m\setminus B(0,R),
			\end{equation}  
			where $2<p_1<2^{\#}$ with $p_1\le p_2\le 2(p_1-1)$ and $2^{\#}\coloneqq\frac{2N-1}{N-1}$.
		\end{enumerate}   
\end{definition}

Assumption~\eqref{eq:W0} encodes the strict triangle inequality for the effective surface
tensions $\omega_{ij}$ and rules out the formation of energetically preferred indirect
transitions between phases. This condition is crucial for the emergence of genuine
multiphase interfaces in the sharp-interface limit and is consistent with the structure
of minimizing isoperimetric clusters. The growth and regularity assumptions
\eqref{eq:W1}, \eqref{eq:W2}, and \eqref{eq:W3} ensure the well-posedness of the variational problem and the
compactness properties needed for $\Gamma$-convergence.
	
For $\varepsilon>0$ and $\mathrm{v}\in\mathbb R^m_{>0}$, we study the
\emph{volume-constrained Allen--Cahn system}
\begin{equation}
\label{eq:AC-PDE}
\tag{$AC_{N,m,\varepsilon,\mathrm{v},g}$}
\begin{dcases}
-\varepsilon\Delta_g u + \dfrac1\varepsilon \nabla W(u)=\lambda,
& \text{in } M,\\[0.3em]
\mathcal{V}_g(u)=\mathrm{v},
\end{dcases}
\end{equation}
where $\lambda\in\mathbb R^m$ is a Lagrange multiplier and $\mathcal{V}_g(u)=\displaystyle\int_M u\,\mathrm dv_g\in\mathbb{R}^m_{>0}$ is the (componentwise) volume functional. Let us define the Hilbert manifold
\begin{equation}\label{eq:constraint-manifold}
 W^{1,2}_{\mathrm{v}}(M,\mathbb R^m)
:=
\Big\{
u\in W^{1,2}(M,\mathbb R^m):
\mathcal V_g(u)=\mathrm{v}
\Big\}   
\end{equation}
We also define the associated energy functional
$\mathcal{AC}_\varepsilon\colon W^{1,2}(M,\mathbb R^m)\to\mathbb R$
as
\begin{equation*}
\mathcal{AC}_\varepsilon(u)
:=
\int_M
\Big(
\frac{\varepsilon}{2}|\nabla_g u|^2
+\frac1\varepsilon W(u)
\Big)\,\mathrm dv_g,
\end{equation*}
where
\[
|\nabla_g u|^2:=\sum_{i=1}^m g(\nabla_g u^i,\nabla_g u^i).
\]
Notice that critical points of $\mathcal{AC}_\varepsilon$ restricted to
$W^{1,2}_{\mathrm{v}}(M,\mathbb R^m)$
are weak solutions of~\eqref{eq:AC-PDE}.
The admissible variations are
\[
TW^{1,2}_{\mathrm{v}}(M,\mathbb R^m)
=
\Big\{
\xi\in W^{1,2}(M,\mathbb R^m):
\int_M \xi\,\mathrm dv_g=0
\Big\}.
\]

The presence of a boundary further enriches this picture, as boundary conditions select
different geometric limits for the interfaces.
When $M$ has a nonempty boundary, the choice of boundary conditions plays a geometric role
in the limiting partition problem. Dirichlet conditions enforce the presence of a fixed
phase along the boundary, while Neumann conditions correspond to free interfaces meeting
$\partial M$ orthogonally in the sharp-interface limit.

In this fashion, it also makes sense to study the Allen--Cahn system in manifolds with boundary.
More precisely, we consider the boundary value problems
\begin{equation}
\label{eq:AC-Dirichlet}
\tag{$AC^{\mathrm{D}}_{N,m,\varepsilon,\mathrm v,g}$}
\begin{cases}
-\varepsilon \Delta u + \dfrac{1}{\varepsilon} \nabla W(u) = \lambda,
& \text{in }  M, \\
\mathcal{V}_g(u) = \mathrm{v}, \\[0.6em]
u = 0,
& \text{on } \partial{M},
\end{cases}
\end{equation}
or
\begin{equation}
\label{eq:AC-Neumann}
\tag{$AC^{\mathrm{N}}_{N,m, \varepsilon,\mathrm v,g}$}
\begin{cases}
-\varepsilon \Delta u + \dfrac{1}{\varepsilon} \nabla W(u) = \lambda,
& \text{in } {M}, \\
\mathcal{V}_g(u) = \mathrm{v}, \\[0.6em]
\dfrac{\partial u}{\partial \nu} = 0,
& \text{on } \partial{M},
\end{cases}
\end{equation}
where $u:M\to\mathbb R^m$ and $\lambda\in\mathbb R^m$ is a Lagrange multiplier.

We define
\[
W^{1,2}_{\mathrm D}(M,\mathbb R^m):=W^{1,2}_0(M,\mathbb R^m),\quad {\rm and} \quad W^{1,2}_{\mathrm N}(M,\mathbb R^m):=W^{1,2}(M,\mathbb R^m).
\]
Also, we set the volume-constraint manifolds as
\begin{equation}\label{eq:constraint-manifold-Dirichlet}
 W^{1,2}_{\mathrm v,\mathrm{D}}(M,\mathbb R^m)
:=\{u\in W^{1,2}_{\mathrm{D}}(M,\mathbb R^m): \mathcal V_g(u)=\mathrm v\}   
\end{equation}
and
\begin{equation}\label{eq:constraint-manifold-Neumann}
 W^{1,2}_{\mathrm v,\mathrm{N}}(M,\mathbb R^m):=\{u\in W^{1,2}_{\mathrm{N}}(M,\mathbb R^m): \mathcal V_g(u)=\mathrm v\}.   
\end{equation}
As before, critical points of $\mathcal{AC}_\varepsilon$ restricted to
$W^{1,2}_{\mathrm v,\mathrm{D}}(M,\mathbb R^m)$ and $W^{1,2}_{\mathrm v,\mathrm{N}}(M,\mathbb R^m)$ are weak solutions to~\eqref{eq:AC-Dirichlet} and \eqref{eq:AC-Neumann}, respectively.

The presence of a boundary introduces additional geometric features that are invisible
in the closed-manifold setting. When $(M,\partial M,g)$ has a nonempty boundary, the choice of boundary
conditions for the Allen--Cahn equation leads to different limiting geometries in the
singular regime. Neumann boundary conditions correspond, at the level of the
$\Gamma$-limit, to free-boundary isoperimetric problems, where interfaces may meet
$\partial M$ orthogonally. Dirichlet boundary conditions, on the other hand, impose
a fixed phase along $\partial M$ and give rise to anchored or constrained interfaces.
These distinctions play a crucial role in multiplicity results and in the structure of
solutions produced by variational and topological methods, and they persist in both the
scalar and vectorial settings.

%====================================================
\subsection{Main results at a glance.}
%====================================================
We present a holistic view of the existing literature on multiplicity results for Allen--Cahn systems.
We prefer to provide a more intuitional rather than a technical description of the results, which we elaborate on in the next sections.
More precisely, in \S~\ref{sec:scalar}, more details about Theorems \ref{thm:intro-scalar-closed}, \ref{thm:intro-scalar-boundary-D}, and \ref{thm:intro-scalar-boundary-N} are provided, and \S~\ref{sec:vectorial} contains further information concerning Theorem \ref{thm:intro-vectorial_bubble} and \ref{thm:intro-vectorial_clusters}.

A guiding principle throughout this survey is that \emph{small-parameter} Allen--Cahn-type problems can be translated, via $\Gamma$-convergence, into \emph{sharp-interface}
isoperimetric questions; the photography method turns the topology of the relevant
configuration space into a multiplicity of critical points.
We organize the main results according to the number of equations and the presence or absence of boundary effects.
More precisely, the number of solutions to volume-constrained Allen--Cahn systems can be
estimated from below in terms of some topological invariants of $M$ and, when present, of its
boundary $\partial M$.

In the scalar theory ($m=1$), the sharp-interface limit is a (relative) isoperimetric problem,
yielding strong geometric rigidity and diameter control; this is the key input that makes
barycenter map constructions and Lusternik--Schnirelmann and Morse theory work effectively.
The first result deals with closed manifolds, {\it i.e., $\partial M=\emptyset$}, and relates the number of solutions to the volume-constrained Allen--Cahn equation with the topology of $M$.

\begin{theoremletter}[{cf.~\cite{MR4396580,MR4644903,MR4314216}}]
\label{thm:intro-scalar-closed}
Let $N,m\in \mathbb{Z}_{>0}$ with $N\ge 2$ and $(M,g)$ be a smooth compact $N$-dimensional Riemannian manifold without boundary.
For sufficiently small $0<\varepsilon,\mathrm{v}\ll1$, Eq. \eqref{eq:AC-PDE} with $m=1$ admits at least $\cat(M)+1$ solutions, and generically, $2\mathcal{P}_1(M)-1$.
\end{theoremletter}

The next results address the multiplicity of solutions to the Allen--Cahn equation when $\partial M\neq\emptyset$. 
In the constrained scalar theory, the presence of $\partial M$ produces a sharp dichotomy. For the Dirichlet boundary value problem, the $\Gamma$-limit analysis remains essentially \emph{interior}, and thus multiplicity detects $\cat(M)$.
\begin{theoremletter}[{cf.~\cite{MR4939669}}]
\label{thm:intro-scalar-boundary-D}
Let $N,m\in \mathbb{Z}_{>0}$ with $N\ge 2$ and $(M,\partial M, g)$ be a smooth compact $N$-dimensional Riemannian manifold with boundary.
For sufficiently small $0<\varepsilon,\mathrm{v}\ll1$, Eq. \eqref{eq:AC-Dirichlet} with $m=1$ admits at least $\cat(M)+1$ solutions, and generically, $2\mathcal{P}_1(M)-1$.
\end{theoremletter}

For the Neumann boundary value problem, however, the singular $\Gamma$-limit becomes the \emph{relative} isoperimetric profile with minimizers asymptotically given by small half-balls attached to $\partial M$.
In this case, multiplicity detects $\cat(\partial M)$ instead of 
$\cat(M)$.
\begin{theoremletter}[{cf.~\cite{MR4939669}}]
\label{thm:intro-scalar-boundary-N}
Let $N,m\in \mathbb{Z}_{>0}$ with $N\ge 2$ and $(M,\partial M, g)$ be a smooth compact $N$-dimensional Riemannian manifold with boundary.
For sufficiently small $0<\varepsilon,\mathrm{v}\ll1$, Eq. \eqref{eq:AC-Neumann} with $m=1$ admits at least $\cat(\partial M)+1$ solutions, and generically, $2\mathcal{P}_1(\partial M)-1$.
\end{theoremletter}

In vectorial systems, the singular $\Gamma$-limit becomes the multi-isoperimetric profile, whose minimizers are \emph{isoperimetric clusters}. 
The absence of a general
classification for minimizers of the multi-isoperimetric profile implies that the uniform diameter
control may fail for the multi-bubble with three or more chambers, producing a genuine obstruction to extending the scalar strategy.
A sharp dichotomy between the double-bubble and multi-bubble (cluster) case emerges.
The double-bubble regime ($m=2$) often retains enough
rigidity to recover multiplicity, which can be stated as follows.
\begin{theoremletter}[{cf.~\cite{MR4701348}}]
\label{thm:intro-vectorial_bubble}
Let $N,m\in \mathbb{Z}_{>0}$ with $N\ge 2$ and $(M,g)$ be a smooth compact parallelizable $N$-dimensional Riemannian manifold without boundary.
For sufficiently small $0<\varepsilon,|\mathrm{v}|\ll1$, Syst. \eqref{eq:AC-PDE} with $m=2$ admits at least $\cat(M)+1$ solutions, and generically, $2\mathcal{P}_1(M)+1$.
\end{theoremletter}

For the multi-bubble regime ($m\ge 3$), the situation is
substantially more delicate.
In fact, the lack of a general classification of minimizers for the multi-isoperimetric profile leads to serious obstructions, and only partial results are currently available under restrictive assumptions.
\begin{theoremletter}[cf.~\cite{andrade2024clusters}]
\label{thm:intro-vectorial_clusters}
Let $N,m\in \mathbb{Z}_{>0}$ with $N\ge 2$ and $(M,g)$ be a smooth compact parallelizable $N$-dimensional Riemannian manifold without boundary. For sufficiently small $0<\varepsilon,\alpha\ll1$, Syst. \eqref{eq:AC-PDE} with $m\geq 3$ and $\bar{\mathrm{v}}=\alpha\mathrm{v}$ admits at least $\cat(M)+1$ solutions, and generically, $2\mathcal{P}_1(M)+1$.
\end{theoremletter}

Throughout this manuscript, when we say that a property holds \emph{generically},
we mean in the sense of \emph{Baire category} with respect to the pair
$(\varepsilon,g)$, where $\varepsilon>0$ is the relaxation parameter and
$g$ is a smooth Riemannian metric on $M$.
More precisely, we consider the product space
\((0,\varepsilon_0)\times \mathrm{Met}^k(M)\),
where $\mathrm{Met}^k(M)$ denotes the space of all $C^k$ Riemannian metrics on $M$
(for $k\gg1$ sufficiently large), endowed with the $C^k$ topology.
A property is said to hold \emph{generically} if it holds on a residual
subset of this space, that is, on a countable intersection of open dense sets.

In this sense, the genericity statements appearing in Theorems
\ref{thm:intro-scalar-closed}–\ref{thm:intro-vectorial_clusters}
mean that, for a residual set of pairs $(\varepsilon,g)$,
\emph{all critical points of the constrained Allen--Cahn functional}
$\mathcal{AC}_\varepsilon$ are \emph{nondegenerate}
in the Morse-theoretical sense.
Equivalently, the second variation of $\mathcal{AC}_\varepsilon$
restricted to the volume constraint manifold
has a trivial kernel at every critical point.

This generic nondegeneracy follows from an infinite-dimensional
transversality argument of Sard--Smale type applied to the Euler--Lagrange
operator associated with the constrained Allen--Cahn system.
A crucial ingredient in this argument is the control of the linearized
operator at critical points, which requires suitable growth bounds on
the Hessian of the potential, a consequence of \eqref{eq:W2}.
For more details on this transversality map construction, see \cite{MR2982783,MR2160744,MR2560131}.
\medskip

%====================================================
\subsection{Organization of the paper}
%====================================================
Section~\ref{sec:sharp-interface} collects some notation and tools from geometric measure theory, which underlie the topological method.
Section~\ref{sec:photography} presents Lusternik--Schnirelmann--Morse theory and the abstract photography construction.
Section~\ref{sec:scalar} treats scalar equations, emphasizing the contrast between the cases with and without boundary.
Section~\ref{sec:vectorial} discusses vectorial systems, which divide into the double-bubble and multi-bubble cases.
Section~\ref{sec:open} lists open problems and perspectives.

%====================================================
%====================================================
\section{Preliminaries}
\label{sec:sharp-interface}
%====================================================
%====================================================
In this section, we discuss basic tools used to construct the photography map.
For convenience, we list the notation that will be used frequently.   
	\begin{itemize}
		\item[] $h=\mathrm{o}(f)$ as $x\rightarrow x_0$ for $x_0\in\mathbb{R}\cup\{\pm\infty\}$, if $\lim_{x\rightarrow x_0}(h/f)(x)=0$;
        \item[] $h=\mathrm{O}(f)$ as $x\rightarrow x_0$ for $x_0\in\mathbb{R}\cup\{\pm\infty\}$, if $-\infty<\lim_{x\rightarrow x_0}(h/f)(x)<\infty$;
        \item[] $\mathbb{Z}_{>0}$ ($\mathbb{Z}_{>0}$) is the set of  positive (nonnegative) integer numbers;
        \item[] $\mathbb{R}_{>0}$ ($\mathbb{R}_{\ge 0}$) is the set of  positive (nonnegative) real numbers;
		\item[] $\mathbb{R}^m_{>0}$ ($\mathbb{R}^m_{\ge 0}$) is the space of vectors with positive (nonnegative) components with $m\in\mathbb{Z}_{>0}$;
		\item[] $(\mathbb{R}^N,\delta)$ is the Euclidean space with its canonical flat metric with $N\ge 2$;
        \item[] $(M,g)$ is a smooth compact $N$-dimensional Riemannian manifold with $N\ge 2$;
		\item[] $B_g(p,r)$ is the geodesic ball with $p\in M$ and $r>0$; we denote $B_{\delta}(p,r)=B(p,r)$;
		\item[] $\Delta_g$ is the Laplace--Beltrami operator;
        \item[] $S_g$ is the scalar curvature.
	\end{itemize}

%====================================================
\subsection{Isoperimetric sets}\label{subsec:GMT}
%====================================================
The singular limit $\varepsilon\to0$ is naturally described in the language of geometric measure theory.
Sets of finite perimeter, reduced boundaries, and varifolds provide the appropriate framework to identify
the limiting interfaces and relate the Allen--Cahn energy to perimeter functionals.
In particular, the geometry of small-volume isoperimetric regions plays a crucial role in the construction
of localized approximate solutions used in the photography method.

For a Borel set $\Omega\subset M$, we denote by $\mathrm{vol}_g(\Omega)$ its Riemannian volume,
by $\chi_\Omega$ its characteristic function, by $\partial^*\Omega$ its reduced boundary, and by
$\mathcal H^{N-1}_g$ the $(N-1)$-dimensional Hausdorff measure associated with $g$.
In this subsection, we introduce the notion of clusters on a Riemannian manifold and formulate the associated
multi-isoperimetric problem. These concepts extend the classical isoperimetric problem ($m=1$) to partitions
with several chambers ($m\ge 2$) and will be used later to describe the $\Gamma$-limit of Allen--Cahn systems.

\begin{definition}[Weighted clusters]\label{def:clusters}
 Let $\mathcal{C}_g^m(M)$ be the class of $m$-finite perimeter subsets of $M$ without intersections, namely
	\begin{equation}
		\label{eq:def-CgmM}
		\mathcal{C}_g^m(M)
		\coloneqq \left\{
		\Omega = (\Omega^1,\dots,\Omega^{m}):
		\begin{aligned}
			\Omega^i\subset M  \text{ is open },
			\, 
			& \Omega^i \cap \Omega^j = \emptyset, \text{ and } \mathcal{H}^{N-1}_g(\partial^* \Omega^i)< \infty
			\\
			& 
			\forall i,j=1,\dots,m \text{ with } i \ne j.
		\end{aligned}
		\right\}.
	\end{equation}
	Each element of $\mathcal{C}_g^m(M)$ will be called 
	\emph{$m$-cluster}
	(or, more simply, \emph{cluster}) and 
	each component of a cluster will be called \emph{a chamber}.
	For every $\Omega\in \mathcal{C}_g^m(M)$,
	we denote by $\Omega^0={\rm int}(M \setminus \bigcup_{i = 1}^m \Omega^i)$ its \emph{exterior} chamber.
	We say that a point $p\in M$
	belongs to a cluster $\Omega \in \mathcal{C}^m_g(M)$,
	and we write $p \in \Omega$,
	if it is in one of its interior chambers,
	hence if $p \in \Omega^i$ for some $i = 1,\dots,m$.
	We also define the
	\emph{diameter} of a cluster as the 
	supremum of the distances between two of its points, namely
	\begin{equation}
		\label{eq:def-diam-cluster}
		\diam(\Omega)\coloneqq
		\sup\Big\{
		\dist_g(p,q):
		p,q \in \Omega
		\Big\}.
	\end{equation}   
\end{definition}

	We define the symmetric difference between two 
	clusters as the cluster obtained by the symmetric 
	differences of their chambers:
	\[
	\Omega_1\triangle\Omega_2
	\coloneqq
	\big(\Omega_1^1\,\triangle\,\Omega_2^1,\,
	\dots,
	\Omega_1^m\,\triangle\,\Omega_2^m\big)\,
	\in \mathcal{C}_g^m(M),
	\qquad \forall \Omega_1,\Omega_2 \in \mathcal{C}^m_g(M).
	\]
	Consequently,
	we define the distance function
	$\dist\colon \mathcal{C}^m_g(M)\times \mathcal{C}^m_g(M) \to \mathbb R_{\geq0}$
	as follows:
	\begin{equation}
		\label{eq:dist-on-clusters}
		\dist(\Omega_1,\Omega_2)
		= \sum_{i = 1}^m
		\mathrm{vol}_g(\Omega_1^i \triangle \Omega_2^i).
	\end{equation}

	For our purposes, we need to associate 
	to every $m$-cluster
	the $\mathbb{R}^m$ valued function that is 
	equal to $\textbf{z}_i$ on the $i$-th chamber,
	and $0$ otherwise, where we recall 
	that $\{\mathbf{z}_1,\dots,\mathbf{z}_m\}$
	are the non-trivial zeros of the potential function $W$.
	Hence, for every $\Omega\in \mathcal{C}^m_g(M)$
	we define $Z_\Omega \in L^1(M,\mathbb{R}^m)$ as 
	\begin{equation}
		\label{eq:def-X_Omega}
		Z_\Omega(p) \coloneqq \sum_{i= 1}^m \mathbf{z}_{i}\chi_{\Omega^i}(p).
	\end{equation}
	Let us notice that $W\circ Z_\Omega$ vanishes 
	a.e. on the whole manifold,
	but $Z_\Omega$ is not a Sobolev function,
	due to the step discontinuities that occur 
	in the interfaces between the chambers
	(excluding, of course, the trivial cases, such as
	$\Omega = \emptyset$).
	
	We define the map $\mathcal{V}_g\colon \mathcal{C}_g^m(M) \to \mathbb{R}^m$,
	which, loosely speaking, is the vectorial volume of clusters,
	as follows:
	\[
	\mathcal{V}_g(\Omega)\coloneqq
	\left(
	\int_{M}
	\chi_{\Omega^1}\de v_{g},
	\dots,
	\int_{M}
	\chi_{\Omega^m}\de v_{g}
	\right)
	\in \mathbb{R}^{m}_{>0}.
	\]
	For every $\mathrm{v}\in \mathbb{R}^m_{>0}$,
	we denote by $\mathcal{C}^m_{g,\mathrm{v}}(M)$
	the subset of $\mathcal{C}^m_g(M)$
	whose elements have vectorial volume $\mathrm{v}\in\mathbb{R}^m_{>0}$,
	namely
	\begin{equation}
		\label{eq:def-CgmMv}
		\mathcal{C}^m_{g,\mathrm{v}}(M)
		\coloneqq
		\Big\{
		\Omega \in \mathcal{C}_g^m(M):
		\int_{M} \chi_{\Omega^i}\de v_{g} = \mathrm{v}^i, \;
		\forall i = 1,\dots, m
		\Big\}.
	\end{equation}
	Let $\omega$ denote the $(m+1)$-dimensional symmetric matrix 
	whose coefficients $\omega_{ij}$ are given by~\eqref{eq:def-omega-ij}
	and let us define 
	the $\omega$-weighted multi-perimeter functional
	$\mathcal{P}^{\omega}_g\colon \mathcal{C}_g^m(M) \to \mathbb{R}$
	as follows:
	\begin{equation}
		\label{eq:def-Pomega}
		\mathcal{P}^\omega_g(\Omega)
		\coloneqq
		\frac{1}{2}
		\sum_{i,j = 0}^{m}
		\omega_{ij}
		\mathcal{H}^{N-1}_g(\partial^*\Omega^{i}\cap \partial^*\Omega^j).
	\end{equation}
	The \emph{weighted multi-isoperimetric problem} asks to minimize
	$\mathcal{P}^\omega_{g}$ among all the 
	clusters with a specified vectorial volume.
	The $\omega$-weighted multi-isoperimetric function
	$\mathcal{I}^{\omega}_g\colon (0,{\rm vol}_g(M))^{m} \to \mathbb{R}$
	associates at each vectorial volume this minimal value,
	hence it is
	\begin{equation}
		\label{eq:def-Iomega}
		\mathcal{I}^{\omega}_{g}(\mathrm{v})
		\coloneqq
		\inf\big\{
		\mathcal{P}^{\omega}_g(\Omega):
		\Omega \in \mathcal{C}^m_{g,\mathrm{v}}(M)
		\big\}.
	\end{equation}
	
		We emphasize that the infimum $\mathcal{I}^{\omega}_{g}(\mathrm{v})$ is always attained,
		and the $m$-clusters that achieve it are known as {\it isoperimetric $\omega$-weighted} $m$-clusters.
		Almgren gave this existence result for the non-weighted case in his classical work~\cite[Chapter VI]{MR420406}: his proof is based on a clever application of the direct method of calculus of variations, hence by considering a minimizing sequence for the multi-perimeter functional, proving that this sequence satisfies the compactness criterion for clusters
		(see~\cite[Proposition 29.5]{MR2976521} for a modern exposition of this compactness result)
		and then relying on the lower semicontinuity of the perimeter functional.
		When different weights $\omega$ are considered, provided they satisfy the immiscibility condition given by~\eqref{eq:W0}, the results in~\cite{MR1070482} ensure that $\mathcal{P}^\omega_g$ is still lower semicontinuous and so $\mathcal{I}^{\omega}_g(\mathrm{v})$ is indeed attained.
		Moreover, 
		in the Euclidean setting, the 
		isoperimetric clusters are regular, and this can be proved by the infiltration lemma and density estimates developed by Leonardi~\cite{leonardi2001}, together with classical results of Almgren and Taylor \cite{MR420406, MR428181}; see in particular \cite[Theorem 3.1]{leonardi2001} and \cite[Theorem 30.1]{MR2976521} for the non-weighted case. Since these arguments are local in nature, they apply verbatim in normal coordinates on a smooth Riemannian manifold, yielding $\mathcal{H}^{N-1}$-a.e. smoothness of the interfaces of weighted minimizing clusters.
	
	All the above definitions, namely
	$\mathcal{C}^m_g(M)$, $\mathcal{P}^{\omega}_g$ and $\mathcal{I}^{\omega}_g$,
	can also be given in the $N$-dimensional
	Euclidean setting, 
	replacing $M$ with $\mathbb{R}^N$.
	The corresponding notions will be denoted 
	replacing $g$ with $\delta$,
	so they are 
	$\mathcal{C}^m_{\delta}$, $\mathcal{P}^{\omega}_{\delta}$
	and $\mathcal{I}^{\omega}_{\delta}$, respectively. To better present our construction, we introduce the notion of almost isoperimetric clusters.
\begin{definition}[Almost minimizing clusters]\label{def:almost-isoperimetric-clusters}
A sequence $(\Omega_k)_{k\in\mathbb N}\subset \mathcal C_g^m(M)$ is called \emph{almost isoperimetric} if
\begin{equation}\label{eq:def-almost-isoperimetric-clusters}
\lim_{k\to\infty}
\frac{\mathcal P_g(\Omega_k)}{\mathcal I_g(\mathrm v_k)} = 1,
\end{equation}
where $\mathrm v_k:=\mathcal V_g(\Omega_k)\in\mathbb R^m_{>0}$.
\end{definition}

If $|\mathrm v_k|\to0$, then $\mathcal I_g(\mathrm v_k)\to0$ (for instance, by taking disjoint unions of small geodesic balls
realizing the prescribed chamber volumes). Hence \eqref{eq:def-almost-isoperimetric-clusters} implies the asymptotic expansion
\begin{equation}\label{eq:almost-isoperimetric-first-term}
\mathcal P_g(\Omega_k)
=
\mathcal I_g(|\mathrm v_k|)
+ o_k\big(\mathcal I_g(|\mathrm v_k|)\big) \quad {\rm as} \quad |\mathrm v_k|\to 0.
\end{equation}

Next, we discuss concentration properties for a multi-isoperimetric cluster enclosing a small volume.
In the double-bubble ($m=2$) regime, one can often localize almost minimizers without losing perimeter at first order,
a key step toward barycenter constructions; see, for instance, \cite[Proposition~12]{MR4701348}.
We will use this type of selection principle later, when discussing the sharp-interface limit of the Allen--Cahn system.

\begin{proposition}\label{prop:concentration-almost-minimizers}

There exists $\mu > 0$ such that the following property holds.
For any infinitesimal sequence $(\alpha_k)_{k\in\mathbb N}\subset\mathbb R_{>0}$,
if $(\Omega_k)_{k\in\mathbb N}\subset\mathcal C_g^m(M)$ is an almost isoperimetric sequence
such that $\mathcal V_g(\Omega_k)=\alpha_k\mathrm v$
for every $k\in\mathbb N$, then there exists a sequence of points $(p_k)_{k\in\mathbb N}\subset M$ such that the sequence $(\Omega_k)_{k\in\mathbb N}\subset\mathcal C_g^m(M)$ defined as
\[
\widetilde{\Omega}_k
:=
\Omega_k\cap B_g\big(p_k,\mu(\alpha_k\|\mathrm v\|)^{1/n}\big)
\]
satisfies
\[
\lim_{k\to\infty}
\frac{|\mathcal V_g(\Omega_k\triangle \widetilde{\Omega}_k)|}{|\mathcal V_g(\Omega_k)|}=0, \quad \lim_{k\to\infty}
\frac{\mathcal P_g(\widetilde{\Omega}_k)}{\mathcal P_g(\Omega_k)}=1,
\]
and
\[
\lim_{k\to\infty}
\frac{\mathrm{vol}_g(\widetilde{\Omega}_k^i)}{\mathrm{vol}_g(\Omega_k^i)}=1,
\qquad \forall\, i=1,\ldots,m.
\]
\end{proposition}
For results concerning properties of isoperimetric clusters, such as regularity of interfaces, we refer the reader to
\cite{MR420406,leonardi2001,MR2976521,MR4470300,MR4588150,MR1402391}.

%====================================================
\subsection{\texorpdfstring{$\Gamma$}{}-convergence}\label{subsec:gamma_convergence}
%====================================================
A fundamental bridge between Allen--Cahn type phase-field models and geometric variational
problems is provided by $\Gamma$-convergence. In the singular limit $\varepsilon\to0$,
diffuse-interface energies concentrate on hypersurfaces and converge, in a variational
sense, to perimeter-type functionals. In the scalar case, this yields the classical
isoperimetric problem, while in vectorial systems, the limit is given by the
multi-isoperimetric problem for clusters.

De Giorgi and Franzoni \cite{MR448194} introduced the notion of $\Gamma$-convergence in 1975, and subsequently developed it in the influential works of Modica–Mortola \cite{MR473971,MR0445362}, Modica \cite{MR866718}, and Kohn–Sternberg \cite{MR985990}. Unlike classical notions such as pointwise or uniform convergence, this variational convergence is well-suited for sequences of functionals arising in minimization problems, as it ensures the convergence of minimizers and minimum values under appropriate conditions. A detailed modern exposition, including the generalized definition for topological spaces, can be found in \cite[Section 2]{BRAIDES2006101}. The notion of $\Gamma$-convergence plays a fundamental role in the variational approach to the Allen--Cahn equation and, in particular, in the photography method for obtaining multiple solutions. Because the energy functional associated with the Allen--Cahn equation $\Gamma$-converges to the perimeter functional of the limiting interface, thus linking phase transition models to geometric variational problems. Historically, the terminology ``$\Gamma$-convergence" arose as a tribute to De Giorgi (hence initially ``$G$-convergence"); however, the letter ``G" was later replaced by its Greek counterpart ``$\Gamma$”.

\begin{definition}[$\Gamma$-convergence]
Let $(X,d)$ be a metric space and let
$\mathcal E_\varepsilon,\mathcal E_0\colon X\to(-\infty,\infty]$.
We say that $(\mathcal E_\varepsilon)$ \emph{$\Gamma$-converges} to $\mathcal E_0$
as $\varepsilon\to0$ (denote by $\Gamma\text{-}\lim_{\varepsilon\to0}\mathcal E_\varepsilon=\mathcal E_0,$
if the following two properties hold:
\begin{enumerate}
\item[\rm (i)] (\emph{liminf inequality})
For every $u\in X$ and every sequence $u_\varepsilon\to u$ in $X$,
\[
\mathcal E_0(u)\le \liminf_{\varepsilon\to0}\mathcal E_\varepsilon(u_\varepsilon).
\]
\item[\rm (ii)] (\emph{limsup inequality})
For every $u\in X$ there exists a sequence $u_\varepsilon\to u$ in $X$
such that
\[
\mathcal E_0(u)\ge \limsup_{\varepsilon\to0}\mathcal E_\varepsilon(u_\varepsilon).
\]
\end{enumerate}
Any sequence $(u_\varepsilon)$ satisfying {\rm (ii)} is called a
\emph{recovery sequence} for $u$.
\end{definition}

Let $\Omega\subset\mathbb R^N$ be a bounded domain with a Lipschitz boundary, and let $W\in C^2(\mathbb{R},\mathbb{R}_{\ge 0})$ be a double-well potential vanishing at $\pm1$. Consider the Allen--Cahn energy functional $\mathcal{AC}_{\varepsilon,\mathrm v}$ associated with a map $u\colon\Omega\to\mathbb R$. For the purpose of $\Gamma$-convergence, it is convenient to extend its domain to $L^1(\Omega,\mathbb R)$ by setting
\[
\mathcal{AC}_\varepsilon(u)
:=\left\{\begin{array}{ll}
\displaystyle\int_\Omega
\left(
\frac{\varepsilon}{2}|\nabla u|^2+\frac1\varepsilon W(u)
\right)\,\mathrm dx,&\mbox{if }u\in W^{1,2}(\Omega,\mathbb R),\\
\infty,&\mbox{if }u\in L^1(\Omega,\mathbb{R})\setminus W^{1,2}(\Omega,\mathbb R).
\end{array}\right.
\]
The limiting functional $\mathcal{AC}_0\colon L^1(\Omega,\mathbb{R})\to[0,\infty]$ is defined by
\[
\mathcal{AC}_0(u)
:=
\begin{cases}
\sigma\mathcal P_\Omega(E),
& \text{if } u=\chi_E
\text{ for some } E\in\mathcal{C}_g(\Omega),\\[4pt]
\infty, & \text{otherwise},
\end{cases}
\]
where $\sigma=2\int_{-1}^1\sqrt{W(s)}\,ds$. Then, as $\varepsilon\to0$, one has
\begin{equation}\label{MM:gamma-convergence}
\Gamma\text{-}\lim_{\varepsilon\to0}\mathcal{AC}_\varepsilon=\mathcal{AC}_0
\quad\text{in }L^1(\Omega,\mathbb{R}).
\end{equation}
This result was established in the seminal works of Modica and Mortola
\cite{MR473971,MR0445362} and Modica \cite{MR866718}, with the construction of optimal
recovery sequences later clarified by Kohn and Sternberg \cite{MR985990}.
Although originally not phrased in the language of $\Gamma$-convergence, these works
form the cornerstone of the modern variational theory of phase transitions.

The classical Modica-Mortola theorem (equivalently to \eqref{MM:gamma-convergence}) establishes that, in the scalar double-well case on Euclidean space, the Allen-Cahn energy $\Gamma$-converges to the perimeter functional. Baldo \cite{MR1051228} improved this result to multi-well potentials. The extension to Riemannian manifolds was developed more recently. For the scalar equation with a double-well potential on manifolds, $\Gamma$-convergence was proved in \cite[Proposition 3.3 and 4.20]{MR4396580}. The first result for vector-valued systems on manifolds, assuming a triple-well potential, appears in \cite[Proposition 11]{MR4701348}. The general multi-well case for the vectorial Allen-Cahn system on Riemannian manifolds was subsequently established in \cite[Proposition 2.4]{andrade2024clusters}. Due to its generality, we present this last result. 

In a framework involving domains with boundary or Riemannian manifolds, the $\Gamma$-convergence preserves the perimeter structure. Under Dirichlet boundary conditions, the limiting functional involves only the interior perimeter, corresponding to the standard relative perimeter in
$\Omega$. Assuming Neumann boundary conditions, boundary contributions may appear in the
Euler--Lagrange equations, but the $\Gamma$-limit remains the perimeter functional, with
minimizers corresponding to relative isoperimetric sets.
Precise formulations for both Dirichlet and Neumann boundary conditions can be found in
\cite{MR1757697,MR1130601,MR930124,MR1620498}.

Let $W\in C^2(\mathbb R^m,\mathbb R_{\ge0})$ be a multi-well potential with exactly
$m+1$ nondegenerate global minima $\{\mathbf z_0,\ldots,\mathbf z_m\}$. Given $\varepsilon>0$ and $\mathrm v\in\mathbb R^m_{>0}$, consider $\mathcal{AC}_{\varepsilon,\mathrm v}$ the vectorial Allen--Cahn energy with the volume constraint defined by

\[
\mathcal{AC}_{\varepsilon,\mathrm{v}}(u):=\left\{\begin{array}{ll}
\int_M
\left(
\frac{\varepsilon}{2}|\nabla_g u|^2+\frac1\varepsilon W(u)
\right)\,\mathrm dv_g,&\mbox{if }u\in W^{1,2}_{\mathrm v}(M,\mathbb R^m)\\
     \infty,&\mbox{if }u\in L^1_{\mathrm{v}}(M,\mathbb R^m)\setminus W^{1,2}_{\mathrm v}(M,\mathbb R^m). 
\end{array}\right.
\]
where $L^1_{\mathrm{v}}(M,\mathbb R^m)=\{u\in L^1(M,\mathbb R^m)\colon \mathcal V_g(u)=\mathrm v\}$. As $\varepsilon\to0$, the functionals $(\mathcal{AC}_{\varepsilon,\mathrm v})$
$\Gamma$-converges in $L^1$ to $\mathcal{AC}_{0,\mathrm v}$ the $\omega$-weighted multi-perimeter functional defined on $\{u\in L^1(M,\mathbb R^m)\colon \mathcal V_g(u)=\mathrm v\}$ in the following way
\begin{equation*}
\mathcal{AC}_{0,\mathrm v}(u):=\left\{\begin{array}{ll}
\mathcal P^\omega_g(\Omega),&\mbox{if }u=Z_\Omega\mbox{ for some }\Omega\in\mathcal C^m_{g,\mathrm v}(M),\\
     \infty,&\mbox{otherwise}. 
\end{array}\right.
\end{equation*}
Here $Z_\Omega$, $\mathcal C^m_{g,\mathrm v}(M)$, and $\mathcal P^\omega_g(\Omega)$ are given in \eqref{eq:def-X_Omega}, \eqref{eq:def-CgmMv}, and \eqref{eq:def-Pomega}, respectively. Moreover, sequences with uniformly bounded energy converge, up to subsequences, to $Z_\Omega$ for some cluster $\Omega\in \mathcal C^{m}_{g,\mathrm v}(M)$. For completeness, we state the precise convergence result established in \cite[Proposition 2.4]{andrade2024clusters}.
	\begin{proposition}
		\label{prop:gamma-convergence}
		For any $\mathrm{v} \in \mathbb{R}^m_{> 0}$
		such that $\sum_{i = 1}^m \mathrm{v}^i < \mathrm{vol}_g(M)$,
		the sequence of functionals
		$(\mathcal{AC}_{\varepsilon,\mathrm{v}})_{\varepsilon > 0}$ 
		$\Gamma$-converges to the weighted perimeter functional
		$\mathcal{P}^{\omega}_{g}\colon \mathcal{C}^m_{g,\mathrm{v}}(M) \to \mathbb{R}$,
		meaning that the following properties hold:
		\begin{itemize}
			\item[{\rm (i)}]
			$($\emph{liminf property}$)$ for any $\Omega \in \mathcal{C}^m_{g,\mathrm{v}}(M)$
			and family
			$(u_{\varepsilon})_{\varepsilon>0}\subset W^{1,2}_\mathrm{v}(M,\mathbb{R}^m)$
			that converges to $Z_\Omega \in L^1(M,\mathbb{R}^m)$,
			namely
			$\lim_{\varepsilon \to 0}
			\norm{u_\varepsilon - Z_\Omega}_{L^1(M,\mathbb{R}^m)}=0$,
			we have
			\begin{equation}
				\label{eq:liminf-prop}
				\liminf_{\varepsilon \to 0}
				\mathcal{AC}_{\varepsilon,\mathrm{v}}(u_\varepsilon)
				\ge \mathcal{P}^\omega_{g}(\Omega);
			\end{equation}
			\item[{\rm (ii)}]
			$($\emph{limsup property}$)$ for any $\Omega\in \mathcal{C}^m_{g,\mathrm{v}}(M)$
			there exists a sequence
			$(u_{\varepsilon})_{\varepsilon>0}\subset W^{1,2}_{\mathrm{v}}(M,\mathbb{R}^m)$
			that converges in $L^1(M,\mathbb{R}^m)$
			to $Z_\Omega$ and 
			\begin{equation}
				\label{eq:limsup-prop}
				\limsup_{\varepsilon \to 0}
				\mathcal{AC}_{\varepsilon,\mathrm{v}}(u_\varepsilon)
				\le \mathcal{P}^\omega_{g}(\Omega).
			\end{equation}
		\end{itemize}
		Moreover, the functionals  
		$\mathcal{AC}_{\varepsilon,\mathrm{v}}$
		are \emph{equicoercive}, that is, for any sequence 
		$(u_{\varepsilon})_{\varepsilon>0}\in W^{1,2}_{\mathrm{v}}(M,\mathbb{R}^m)$
		such that
		\begin{equation}\label{prop:equicoerciveness}
			\limsup_{\varepsilon\to 0}
			\mathcal{AC}_{\varepsilon,\mathrm{v}}(u_\varepsilon)
			< \infty,
		\end{equation}
		there exists $\Omega \in \mathcal{C}^m_{g,\mathrm{v}}(M)$
		such that $(u_{\varepsilon})_{\varepsilon>0}$
		converges in $L^1(M,\mathbb R^m)$ to $Z_\Omega$,
		up to subsequences.
	\end{proposition}

%====================================================
%====================================================
\section{Topological method}
%====================================================
%====================================================
\label{sec:photography}
As a guiding example, let us consider the Allen--Cahn functional under a scalar
volume constraint, that is, for a fixed $\mathrm{v} \in \mathbb{R}_{> 0}$, we study the functional 
\( \mathcal{AC}_{\varepsilon,\mathrm{v}}: 
    W^{1,2}_{\mathrm{v}}(M,\mathbb{R}) \to \mathbb{R}.
\)
By a standard application of the direct method,
if $(M,g)$ is a compact manifold, the existence of a minimizer for
Allen-Cahn functional is always guaranteed.
A natural question is therefore the existence of critical points
different from the minima.
As is customary in critical point theory,
one relies on topological arguments to obtain lower bounds on the number of
critical points.
Among the most classical tools are Lusternik-Schnirelmann theory and Morse theory.

%====================================================
\subsection{Lusternik--Schnirelmann--Morse theory}
The \emph{Lusternik-Schnirelmann category} $\mathrm{cat}(X)$ of a topological space $X$
is a homotopy invariant that measures the minimal number of contractible sets
needed to cover $X$, namely
\begin{equation}\label{def:LS_category}
\mathrm{cat}(X)
\coloneqq
\min \Big\{
k \in \mathbb{N}: X = \bigcup_{i=1}^k A_i,
\ \text{with each } A_i \text{ contractible in } X
\Big\}.
\end{equation}
From a variational viewpoint, it provides a natural lower bound on the number of
critical points of a smooth functional whose sublevels reflect the topology of $X$,
since each nontrivial topological feature gives rise to a min-max construction,
and hence to a critical point.
More precisely, if $f \colon X \to \mathbb{R}$ is a $C^1$ functional, then for any
$c \in \mathbb{R}$, the number of critical points in the sublevel
\[
f^c := \{u\in X : f(x) \le c\}
\]
is at least $\mathrm{cat}(f^c)$.
We recall that the inclusion $f^c \subset X$ does not imply
$\mathrm{cat}(f^c) \le \mathrm{cat}(X)$, since the contraction must be 
verified within the sublevel itself.
The Lusternik-Schnirelmann approach dates back to the seminal works of Lusternik and Schnirelmann
in~\cite{LS34},
where the notion of category was introduced in connection with classical
results such as the Borsuk--Ulam Theorem (we refer to~\cite{CLOT03} for a modern exposition).
More recently, it has been successfully employed to provide lower bounds on the
number of periodic orbits of Hamiltonian systems
(cf.~\cite{MR4113326,MR4305588,MR4395723,MR3412384}),
as well as to obtain existence and multiplicity results for critical points of
indefinite Lagrangian action functionals admitting symmetries
(see, e.g.,~\cite{MR4525720,MR4772577}).

A sharper, but somewhat less general, method for obtaining multiplicity results for critical points is provided by the \emph{Morse theory}.
Roughly speaking, this approach relates the topology of a manifold
to the local behavior of suitably regular functions near their critical points.
We limit ourselves here to a very brief and informal presentation,
aimed only at conveying the main guiding principles of the theory.
A proper treatment of Morse theory
would deserve a much more detailed discussion,
to which we refer the interested reader in the standard literature~\cite{MR163331} for a classical reference, and~\cite{MR1849220}
for a modern exposition.

Let $\mathfrak{M}$ be a (finite or infinite-dimensional) Hilbert manifold,
and let $f \colon \mathfrak{M} \to \mathbb{R}$ be a $C^2$ function
with nondegenerate critical points.
Denoting by $K_f \subset \mathfrak{M}$ the set of critical points of $f$,
a point $x \in K_f$ is said to be \textit{nondegenerate} if the kernel of the bilinear form
$\mathrm{d}^2 f(x)[\cdot,\cdot]$ is trivial.
The \emph{Morse index} of $x \in K_f$, denoted by $m_f(x)$,
is the maximal dimension of a subspace of $T_x\mathfrak{M}$ on which
$\mathrm{d}^2 f(x)$ is negative definite.
If all critical points of $f$ are nondegenerate, then $f$ is called a
\emph{Morse function}, and one can associate to it the \emph{Morse polynomial} defined by
\[
    \mathcal{P}_{f}(t)
    =
    \sum_{x \in K_f} t^{m_f(x)}
    =
    \sum_{k = 0}^\infty a_k\, t^k,
\]
where $a_k$ denotes the number of critical points of $f$
with Morse index equal to $k$.

The Morse theory provides a relation between the Morse polynomial of a function and the Poincar\'e polynomial of the manifold, where the last is given by
\[
    \mathcal{P}_{\mathfrak{M}}(t)
    :=
    \sum_{k \ge 0} \beta_k(\mathfrak{M})\, t^k,
\]
where $\beta_k(\mathfrak{M})$ are the Betti numbers of $\mathfrak{M}$.
More precisely, the Morse inequalities imply that
\[
    \mathcal{P}_{f}(t)
    =
    \mathcal{P}_{\mathfrak{M}}(t)
    + (1+t)\,Q(t),
\]
for some polynomial $Q(t)$ with nonnegative integer coefficients,
depending on the functional $f\in C^2(\mathfrak{M},\mathbb{R})$.
In this way, the topology of a manifold,
encoded in its Poincar\'e polynomial,
gives lower bounds on the number of critical points of any Morse function.
Moreover, in contrast with Lusternik--Schnirelmann theory,
the above relation also provides information on the nature of the critical points,
as it takes into account their Morse indices.

%====================================================
\subsection{Photography method}
In view of the considerations above, and returning to our original question of estimating the number of critical points of the Allen--Cahn functional under a volume constraint, the application of Lusternik-Schnirelmann and Morse techniques
require some control on the topology of the functional space,
or at least of an appropriate sublevel.
If $M$ is a closed (compact and without boundary) Riemannian manifold,
it is natural to expect that a minimizer of the Allen--Cahn functional under a small volume constraint
is concentrated in a small region of the manifold.
Indeed, $\Gamma$-convergence shows that, as $\varepsilon \to 0$,
the minimizers $u_\varepsilon$ of $\mathcal{AC}_\varepsilon$
converge to an isoperimetric region with the same volume,
whose existence follows from the direct method
and, in the compact setting, such isoperimetric regions consist of a single connected component.
As a consequence, for $\varepsilon$ sufficiently small,
the same qualitative behavior is expected for the minimizers 
of $\mathcal{AC}_\varepsilon$.

More importantly, this phenomenon is stable in the following sense:
an \emph{almost isoperimetric} region under a small volume constraint
has \emph{almost all} of its mass concentrated in a small ball
(see Proposition~\ref{prop:concentration-almost-minimizers})
and,
via $\Gamma$-convergence, the same applies for the Allen-Cahn functional:
there exists a sublevel
\[
    \mathcal{AC}_\varepsilon^c := \{u : \mathcal{AC}_\varepsilon(u) \le c\}
\]
such that every function $u \in \mathcal{AC}_\varepsilon^c$
is essentially concentrated near a point of the manifold.
More precisely, there exist $p \in \Omega$ (depending on $u$) and a small radius $r > 0$
(independent of $u$) such that
\[
    \frac{\int_{B(p,r)} u \, \mathrm{d}x}
         {\int_{\Omega} u \, \mathrm{d}x}
    \ge 1 - \eta,
\]
where $\eta > 0$ can be taken arbitrarily small, depending on $c$.

This concentration property is the key ingredient for extracting topological information
on the sublevel $\mathcal{AC}_\varepsilon^c$
and, consequently, for obtaining lower bounds on the number of critical points of the functional.
Indeed, if $c$ is sufficiently close to the minimum of $\mathcal{AC}_\varepsilon$,
the above concentration property allows to associate to each function in the sublevel $\mathcal{AC}_\varepsilon^c$
a point of the manifold in a canonical and continuous way.
If, conversely, one can also associate to each point of the manifold
a function in $\mathcal{AC}_\varepsilon^c$,
again in a continuous and canonical manner,
and if the composition of these two maps is homotopic to the identity,
then the topology of the sublevel reflects that of the manifold itself.
As a consequence, one obtains lower bounds on the number of critical points
in terms of topological invariants of the manifold.
In other words, we need to find a real number $c \in \mathbb{R}$
and two continuous maps
\[
\varphi:M \longrightarrow \mathcal{AC}_\varepsilon^c,
\qquad 
\beta: \mathcal{AC}_\varepsilon^c \longrightarrow M, 
\]
such that 
\[
    \beta \circ \varphi:M \longrightarrow M
    \text{ is homotopic to the identity map}.
\]
This strategy is known as the \emph{photography method}:
if the underlying manifold can be continuously embedded into a suitable sublevel of the energy functional,
then the sublevel is topologically no simpler than the manifold itself,
even if it is much larger.
In this sense, the sublevel contains a ``photograph'' of the manifold,
capturing its essential topological features.
More formally, we have the following
the general statement of the abstract photography method.
% and its proof can be found in~\cite[Section 2]{MR4701348}
\begin{theorem}
	\label{theorem:abstract-photography}
	Let $X$ be a topological space,
	$\mathfrak{M}$ be a $C^2$-Hilbert manifold,
	$\mathcal{E}\colon \mathfrak{M}\to\mathbb{R}$ be a $C^1$-functional,
	and $\mathcal{E}^c\coloneqq\{u\in\mathfrak{M} : \mathcal{E}(u)\le c\}$
	be a sublevel set for some $c\in\mathbb R$.
	Assume that the following properties hold:
	\begin{enumerate}
		\item[($E_1$)]
			$\inf_{u\in\mathfrak{M}}\mathcal{E}(u)>-\infty$;
		\item[($E_2$)]
			$\mathcal{E}$ satisfies the Palais-Smale (PS) condition;
		\item[($E_3$)]
			There exist $c\in\mathbb R$ and two continuous maps
			$\Psi_{R}\colon X\rightarrow \mathcal{E}^c$ and
			$\Psi_{L}\colon\mathcal{E}^c\rightarrow X$
			such that $\Psi_{L}\circ\Psi_{R}$ is homotopic
			to the identity map of $X$.
	\end{enumerate} 
	Then, the number of critical points in $\mathcal{E}^c$ 
	is greater than $\mathrm{cat}(X)$,
	and if $\mathfrak{M}$ is contractible and $\mathrm{cat}(X)>1$,
	there exists at least another critical point of $\mathcal{E}$ outside
	$\mathcal{E}^c$.
	Moreover,
	if all the critical points are nondegenerate,
	there exists $c_0\in(c,\infty)$ such that 
	$\mathcal{E}^{c}$ contains $\mathcal{P}_1(X)$ critical points
	and $\mathcal{E}^{c_0}\setminus \mathcal{E}^{c}$
	contains $\mathcal{P}_1(X)-1$ critical points
	if counted with their multiplicity.
	More precisely, the following relation holds:
	\begin{equation}
		\label{eq:morserelation}
		\sum_{u\in {\rm Crit(\mathcal E)}} t^{\mu(u)}
		=\mathcal{P}_t(X)+t[\mathcal{P}_t(X)-1]+(1+t)\mathcal{Q}(t),
	\end{equation}
	where $\mathcal{Q}(t)$ is a polynomial
	with nonnegative integer coefficients,
	${\rm Crit}(\mathcal{E})$
	denotes the set of critical points of $\mathcal{E}$
	and $\mu(u)$ denotes the (numerical) Morse index
	of $u$, {\it i.e.}, the dimension of the maximal subspace
	on which the bilinear form
	$\mathrm{d}^2\mathcal{E}_{\varepsilon}(u)[\cdot,\cdot]$
	is negative-definite.
\end{theorem}

The photography method was introduced in the early 1990s by Benci and
Cerami.
In particular, in~\cite{MR1088278} the authors introduced this approach
to obtain lower bounds on the number of solutions of the following family
of elliptic problems:
\[
\begin{cases}
    -\Delta u + \lambda u = u^{p-1}, & \text{in } \Omega,\\
    u > 0, & \text{in } \Omega,\\
    u = 0, & \text{on } \partial\Omega,
\end{cases}
\]
under suitable assumptions on $p \in (2,2^*)$ and $\lambda \in \mathbb{R}$.
In that work, the estimate was obtained by relying solely on
Lusternik--Schnirelmann theory.
The same authors subsequently combined the photography method with
Morse theory in~\cite{MR1384393}
(see also~\cite{MR1322324}), obtaining sharper multiplicity results
for the same class of problems, allowing for more general nonlinearities
instead of the pure-power one.

Since its introduction, the photography method has attracted considerable
attention and has given rise to a broad literature.
Its applicability relies on the presence of a suitable concentration
property for the sublevels of the energy functional, and, beyond the
original setting of semilinear elliptic problems with variational structure,
it has been applied in a wide range of contexts.
Prominent examples include nonlinear Schr\"odinger-type equations and
Schr\"odinger--Poisson systems, where the method has been used to detect
multiple standing wave solutions and semiclassical states
(cf.~\cite{MR1646619,MR1734531,MR4924234}).
Geometric applications arise in conformal problems, such as the Yamabe equation
and equations involving prescribed $Q$-curvature
(cf.~\cite{MR3912791,MR4761862}).
The same approach has also been employed in the study of nonlocal and coupled
elliptic problems, including Kirchhoff-type and Boussinesq-type equations
(cf.~\cite{MR4862498}).

Because of the natural concentration properties exhibited by low-energy sublevels of the Allen--Cahn functional, this framework has been effectively
implemented to obtain multiplicity results for the Allen--Cahn equation, as we review in the next section.

%====================================================
%====================================================
\section{Scalar equations}\label{sec:scalar}
%====================================================
%====================================================
This section briefly explains the proof of Theorem~\ref{thm:intro-scalar-closed}, which addresses the case of closed manifolds, and the proofs of Theorems~\ref{thm:intro-scalar-boundary-D} and~\ref{thm:intro-scalar-boundary-N}, which apply to manifolds with boundary.

%====================================================
\subsection{Manifolds without boundary}
%====================================================
As previously explained, on closed Riemannian manifolds, the scalar
Allen--Cahn functional under a small volume constraint is expected to exhibit a strong concentration behavior for low-energy configurations, which localize near small geodesic balls.
The photography method exploits this property to control the topology of a sublevel, so that
Lusternik--Schnirelmann and Morse theories then apply.

The first multiplicity result for the Allen--Cahn equation obtained via the
photography method appears in~\cite{MR4073210}.
In this seminal work, the domain $M$ is an open and bounded subset of
$\mathbb{R}^N$, and the double-well potential $W$ is asymmetric, in the sense
that it possesses two non-degenerate local minima with different values.
More precisely, the potential satisfies
\[
W(0)=W'(0)=0,\quad W''(0)>0,
\]
and there exists $s_0>0$ such that
\[
W'(s_0)=0,\quad W''(s_0)>0,\quad W(s_0)<0.
\]
Under these assumptions, the main result concerns multiplicity
for critical points of the Allen--Cahn functional
$\mathcal{AC}_\varepsilon$ under a scalar volume constraint, namely on the set
\(W^{1,2}_{\mathrm v}(M,\mathbb R^m)\) defined as \eqref{eq:constraint-manifold}.

Although the proof is based on the photography method, it differs from
subsequent works because of the asymmetry of the potential $W$.
This feature allows for the introduction of a suitable auxiliary minimization
problem, which plays a central role in the construction.
More precisely, for any fixed $\tau>0$, the authors consider the minimization problem for $\mathcal{AC}_\varepsilon$ over the convex set
\[
K_\tau \coloneqq \Big\{u \in W^{1,2}(\mathbb{R}^N) :
u \ge 0,\ \int_{\mathbb{R}^N} u\, \mathrm{d} x \le \tau \Big\}.
\]
We highlight that this problem is defined on the whole $\mathbb{R}^N$,
and it serves as a model for the subsequent construction.
Indeed, a careful and nontrivial analysis of this auxiliary problem shows that its
minimizers have compact support.
These compactly supported solutions are then suitably rescaled to define a
photography map
\[
    \varphi_{\varepsilon,\mathrm{v}} : M \longrightarrow \mathcal{AC}_\varepsilon^c
    \subset W^{1,2}_{\mathrm v}(\mathbb{R}^N,\mathbb R^m),
\]
where the energy level $c\in\mathbb{R}$ can be estimated in terms of the auxiliary problem.
This approach bypasses the use of $\Gamma$-convergence techniques, which will
instead play a central role in later developments of the theory.
Then, a corresponding barycenter map is defined by
\[
\beta(u) :=
\frac{\int_M x\,|u(x)|\,\mathrm{d}x}
     {\int_M |u(x)|\,\mathrm{d}x},
\]
where the denominator is nonzero due to the volume constraint.
This definition is natural in the present setting, since $M \subset
\mathbb{R}^N$ and no geometric obstruction arises.

The subsequent contribution on this topic is~\cite{MR4396580}
(see also the corrigendum~\cite{MR4644903}).
In that work, $(M,g)$ is a closed (compact and without
boundary) $N$-dimensional Riemannian manifold, the double-well potential $W$ is (subcritical) symmetric
and satisfies suitable structural assumptions
(for simplicity, one can consider a standard symmetric double-well potential $W(u)=\tfrac{1}{4}(1-u^2)^2$),
and the volume constraint is the standard scalar one, so that the functional
space $W^{1,2}_\mathrm{v}(M,\mathbb{R})\subset W^{1,2}(M,\mathbb{R})$ is given by \eqref{eq:constraint-manifold}.

The main result \cite[Theorem~2]{MR4396580} is essentially analogous to Theorem~\ref{thm:BNP-firstWork},
but now set in a Riemannian framework.
In contrast with the previous work, the proof strategy relies heavily on the
$\Gamma$-convergence of the Allen--Cahn functional to the perimeter and on
concentration properties of almost-isoperimetric sets.

We briefly outline the main ingredients of the argument, which will also
underlie the subsequent results discussed below, up to minor adaptations
dictated by the specific setting.
In the small-volume regime, the perimeter of a geodesic ball centered at
$p\in M$ with volume $\mathrm v\in\mathbb{R}_{> 0}$, denoted by $B_g(p,r_\mathrm v(p))$, admits the expansion
\begin{equation}
    \label{eq:perimeter-expansion}
    \mathcal{P}_g\big(B_g(p,r_\mathrm v(p))\big)
    = c_N \mathrm{v}^{\frac{N-1}{N}}
    - \gamma_N S_g(p)\, \mathrm{v}^{\frac{N+1}{N}}
    + o\big(\mathrm{v}^{\frac{N+1}{N}}\big) \quad {\rm as} \quad \mathrm{v}\to 0,
\end{equation}
where $c_N>0$ is the Euclidean isoperimetric constant,
$\gamma_N>0$ depends only on the dimension,
and $S_g(p)$ denotes the scalar curvature of $M$ at $p$
(cf.~\cite[Corollary~2]{MR4130849}).
By compactness, this expansion holds uniformly with respect to
$p\in M$.
By the lim-sup property of $\Gamma$-convergence, for $0<\varepsilon\ll1$
sufficiently small, one can approximate each ball
$B_g(p,r_\mathrm v(p))$ by a Sobolev function
$\varphi_{\varepsilon,\mathrm v}(p)\in W^{1,2}_\mathrm{v}(M,\mathbb{R})$
satisfying
\[
\lim_{\varepsilon\to0}
\mathcal{AC}_\varepsilon\big(\varphi_{\varepsilon,\mathrm{v}}(p)\big)
= \sigma\,\mathcal{P}_g\big(B_g(p,r_\mathrm{v}(p))\big),
\]
where $\sigma>0$ is a constant depending only on the potential $W$.
This provides the construction of a photography map
$\varphi_{\varepsilon,\mathrm{v}}:M\to W^{1,2}_\mathrm{v}(M,\mathbb{R})$.

Moreover, using the compactness of $M$ again, one can estimate a value
$c(\mathrm{v})\in\mathbb{R}_{>0}$ such that
\[
\varphi_{\varepsilon,\mathrm{v}}(M)\subset\mathcal{AC}_\varepsilon^{c(\mathrm{v})}.
\]
Indeed, it suffices to set
\[
c(\mathrm{v})
= \sigma\, c_N \mathrm{v}^{\frac{N-1}{N}}
+ C \mathrm{v}^{\frac{N+1}{N}}
= \sigma \mathcal{I}_{\delta}(\mathrm{v})
+  \mathrm{v}^{\frac{N+1}{N}},
\]
for some constant $C>0$ depending on the scalar curvature
$S_g$.
This estimate is uniform with respect to $\varepsilon$, provided $0<\varepsilon\ll1$
is sufficiently small.
The above bound allows one to define a barycenter map on the sublevel
$\mathcal{AC}_\varepsilon^{c(\mathrm{v})}$.
Indeed, one has
\begin{equation}
    \label{eq:cvtoIv}
    \lim_{\mathrm{v}\to0}\frac{c(\mathrm{v})}{\sigma \mathcal{I}_g(\mathrm{v})}
    = \lim_{\mathrm{v}\to0}\frac{c(\mathrm{v})}{\sigma \mathcal{I}_{\delta}(\mathrm{v})}
    = 1,
\end{equation}
which is nothing but the equivalent of 
~\eqref{eq:def-almost-isoperimetric-clusters} in the definition of almost 
isoperimetric sequence (cf. Definition~\ref{def:almost-isoperimetric-clusters}).
As a consequence, for $0<\mathrm{v}\ll1$ sufficiently small and $\varepsilon\to0$, any family
$u_\varepsilon\in\mathcal{AC}_\varepsilon^{c(\mathrm{v})}$
converges in $L^1$ to the characteristic function of a finite-perimeter set
$\Omega_\mathrm{v}\in\mathcal{C}_g(M)$ with enclosed volume $\mathrm{v}\in\mathbb{R}_{>0}$ satisfying
\[
\lim_{\mathrm{v}\to0}\frac{\mathcal{P}_g(\Omega_\mathrm{v})}{\mathcal{I}_g(\mathrm{v})}=1.
\]
In particular, $\{\Omega_\mathrm{v}\}_{\mathrm{v}>0}$ is an almost-isoperimetric family as $\mathrm{v}\to0$,
and hence each $\Omega_\mathrm{v}$ has almost all its volume concentrated in a small
ball (cf.~\cite[Theorem~4.1]{MR4130849}).
The same concentration property is inherited by the approximating functions
$u_\varepsilon$.
More formally, by a diagonal argument as $\mathrm{v}\to0$ and $\varepsilon\to0$, one obtains the following result.
\begin{lemma}[{cf.~\cite[Lemma~4.3]{MR4644903}}]
\label{lem:concentration}
There exists $\mu>0$ such that for every $\rho\in(0,1)$ there exists $\mathrm{v}^*>0$ such that for every
$\mathrm{v}\in(0,\mathrm{v}^*)$ there exists $\varepsilon^*>0$ with the property that for every
$\varepsilon\in(0,\varepsilon^*)$ and
$u\in\mathcal{AC}_\varepsilon^{c(\mathrm{v})}$
there exists a point $p_u\in M$ such that
\[
\int_{B_g(p_u,\mu \mathrm{v}^{1/N})} u\,\mathrm{d}v_g
\ge (1-\rho)\,\mathrm{v}.
\]
\end{lemma}

In other words, every function in the sublevel
$\mathcal{AC}_\varepsilon^{c(\mathrm{v})}$
has almost all its mass concentrated in a small ball.
Consequently, by Nash's embedding theorem, we may view $M$ as a subset of
$\mathbb{R}^N$ and define the extrinsic barycenter map
\[
\beta^*:\mathcal{AC}_\varepsilon^{c(\mathrm{v})}\to\mathbb{R}^N,
\qquad
\beta^*(u)
:=\frac{\int_M x\,|u(x)|\,\mathrm{d}x}
        {\int_M |u(x)|\,\mathrm{d}x}.
\]
For $0<\mathrm{v},\varepsilon\ll1$ sufficiently small, the image of $\beta^*$ lies in a
tubular neighborhood of $M$, so that a projection yields a barycenter map
\[
\beta:\mathcal{AC}_\varepsilon^{c(\mathrm{v})}\to M.
\]
Since each $\varphi_{\varepsilon,\mathrm{v}}(p)$ is obtained via the recovery sequence
construction of the small ball
$B_g(p,r_\mathrm{v}(p))$, one has
$\beta(\varphi_{\varepsilon,\mathrm{v}}(p))\approx p$.
An explicit homotopy between
$\beta\circ\varphi_{\varepsilon,\mathrm{v}}: M\to M$
and the identity can then be constructed using the exponential map, which
concludes the argument.

Summarizing, the whole construction relies on a few key properties of the
\emph{limit geometric problem}, which make the photography method applicable.
In the Allen--Cahn setting, these properties are inherited by the
$\varepsilon$-dependent functional through $\Gamma$-convergence arguments.
More precisely, the essential ingredients are:
\begin{itemize}
    \item[{\rm (i)}] a sharp small-volume expansion of the perimeter of geodesic balls
    (cf.~\eqref{eq:perimeter-expansion}),
    which provides a natural candidate for isoperimetric sets;
    \item[{\rm (ii)}] the convergence, as $\mathrm{v}\to0$, of this expansion to the isoperimetric
    profile, in the sense of~\eqref{eq:cvtoIv}, which yields a concentration
    property for low-energy sublevels, as in Lemma~\ref{lem:concentration}.
\end{itemize}
With these ingredients, we have
\begin{theorem}[{cf.~\cite[Theorem~1.3]{MR4073210}}]
\label{thm:BNP-firstWork}
Let $N,m\in\mathbb{Z}_{>0}$ with $N\ge2$ and $({M},g)$ be a compact $N$-dimensional Riemannian manifold without boundary.
Then, there exists $\mathrm{v}^*>0$ such that for every $\mathrm{v}\in(0,\mathrm{v}^*)$ there exists
$\varepsilon^*=\varepsilon^*(M,g,\mathrm{v})>0$ satisfying that
for every $\varepsilon\in(0,\varepsilon^*)$, Eq.~\eqref{eq:AC-PDE} with $m=1$ admits at least $\mathrm{cat}({M})+1$ distinct solutions.
Moreover, in the nondegenerate setting, one can find at least $2\mathcal{P}_1({M})-1$ distinct solutions.
\end{theorem}

%====================================================
\subsection{Manifolds with boundary}
%====================================================
Since the above properties also hold in the case of compact Riemannian manifolds
with boundary, analogous multiplicity results can be obtained in this setting
as well, as shown in~\cite{MR4939669}.
In this framework, however, the number of critical points depends on the
boundary conditions imposed, and both Dirichlet and Neumann conditions can be
considered, namely \eqref{eq:AC-Dirichlet} and \eqref{eq:AC-Neumann}.
In what follows, we provide the precise statement and briefly explain the constructions in Theorem~\ref{thm:intro-scalar-boundary-D} and Theorem~\ref{thm:intro-scalar-boundary-N}.

When Dirichlet boundary conditions are prescribed, the Allen--Cahn functional
is naturally defined on a volume-constrained subset
$W^{1,2}_{\mathrm{v},\mathrm{D}}({M},\mathbb{R})\subset W^{1,2}_0({M},\mathbb{R})$ defined as \eqref{eq:constraint-manifold-Dirichlet}. In this case, $(\mathcal{AC}_\varepsilon)$ $\Gamma$-converges to the standard perimeter functional, and the arguments developed in~\cite{MR4396580} can be applied with only minor technical modifications to account for the presence of the boundary. As a consequence, the resulting multiplicity estimates are again expressed in
terms of topological invariants of $M$.
\begin{theorem}[{Dirichlet case, cf.~\cite[Theorem B]{MR4939669}}]
\label{thm:dirichlet-summary}
Let $N,m\in\mathbb{Z}_{>0}$ with $N\ge2$ and $({M},g)$ be a compact $N$-dimensional Riemannian manifold with boundary.
Then, there exists $\mathrm{v}^*>0$ such that for every $v\in(0,\mathrm{v}^*)$ there exists
$\varepsilon^*=\varepsilon^*(M,g,\mathrm{v})>0$ satisfying that
for every $\varepsilon\in(0,\varepsilon^*)$, Eq.~\eqref{eq:AC-Dirichlet} with $m=1$ admits at least $\mathrm{cat}({M})+1$ distinct solutions.
Moreover, in the nondegenerate setting, one can find at least $2\mathcal{P}_1({M})-1$ distinct solutions.
\end{theorem}

A different behavior arises in the case of Neumann boundary conditions,
namely when $\partial u/\partial\nu=0$ on $\partial{M}$.
In this setting, the $\Gamma$-limit is given by the \emph{relative} perimeter,
which only accounts for the interior part of the reduced boundary, that is, for any finite-perimeter set $\Omega\in \mathcal{C}^1_g(M)$, one has
\[
\mathcal{P}^*_g(\Omega)
= \mathcal{H}^{N-1}_g\big(\partial^*\Omega \cap \mathrm{int}({M})\big).
\]
Accordingly, the natural candidates for isoperimetric sets are
\emph{half-bubbles} attached to the boundary of the manifold, and the
photography map is defined on $\partial{M}$, namely
\[
\varphi_{\varepsilon,\mathrm{v}} : \partial{M} \longrightarrow W^{1,2}_{\mathrm{v},\mathrm{N}}(M,\mathbb{R}),
\]
where $W^{1,2}_{\mathrm{v},\mathrm{N}}(M,\mathbb{R})\subset W^{1,2}(M,\mathbb{R})$ is defined as \eqref{eq:constraint-manifold-Neumann}.

Sharp small-volume expansions for the perimeter of such half-bubbles are
available in~\cite{MR2587431}, and are analogous to~\eqref{eq:perimeter-expansion},
with the leading term given by the isoperimetric profile of the Euclidean
half-space and the second-order term involving the mean curvature of
$\partial{M}$.
Moreover, compactness results from~\cite{MR4467099} yield a concentration
property for the relevant sublevels of the energy functional, showing that
the volume concentrates near the boundary.
As a consequence, the final multiplicity result is expressed in terms of
topological invariants of the boundary.
\begin{theorem}[{Neumann case, cf.~\cite[Theorem A]{MR4939669}}]
\label{thm:Neumann-summary}
Let $N,m\in\mathbb{Z}_{>0}$ with $N\ge2$ and $({M},g)$ be a compact $N$-dimensional Riemannian manifold with boundary.
Then, there exists $\mathrm{v}^*>0$ such that for every $v\in(0,\mathrm{v}^*)$ there exists
$\varepsilon^*=\varepsilon^*(M,g,\mathrm{v})>0$ satisfying that
for every $\varepsilon\in(0,\varepsilon^*)$, Eq.~\eqref{eq:AC-Neumann} with $m=1$ admits at least $\mathrm{cat}(\partial{M})+1$ distinct solutions.
Moreover, in the nondegenerate setting, one can find at least $2\mathcal{P}_1(\partial{M})-1$ distinct solutions.
\end{theorem}

%====================================================
%====================================================
\section{Vectorial systems}\label{sec:vectorial}
%====================================================
%====================================================
Allen--Cahn systems correspond, in the singular limit, to
multi-phase partition problems, namely, to perimeter minimization problems for
clusters with prescribed volumes.
As in the scalar case, minimizers of the limiting problem exist and enjoy
basic localization properties in the small-volume regime.
However, the geometry of minimizing clusters is in general far more delicate:
unlike isoperimetric regions, which are well understood and essentially rigid
for small volumes, minimizers for multi-chamber clusters need not admit an
explicit description, nor a canonical choice depending continuously on the
base point.

From the viewpoint of the photography method, this lack of rigidity becomes
relevant when one attempts to construct a global family of model clusters and
to enforce the prescribed volumes after projecting them onto the manifold.
In the multi-phase case, these steps typically require additional tools,
such as volume-fixing procedures and infiltration-type lemmas controlling the
perimeter cost of local modifications.

%====================================================
\subsection{Double-bubble case}
%====================================================
When only two phases are prescribed (with the other playing the role of
the exterior region), the limiting isoperimetric problem reduces to the
classical double-bubble problem in Euclidean space, for which minimizers are
explicitly classified.
This geometric rigidity makes it possible to implement the photography method
in a manner closer to the scalar setting.
Subsequently, we provide the precise statement and briefly explain the constructions in Theorem~\ref{thm:intro-vectorial_bubble}.

In~\cite{MR4701348}, a multiplicity result of the same flavor is established
for the Allen--Cahn functional in the vectorial setting, namely for maps
$u \in W^{1,2}({M},\mathbb{R}^2)$, where ${M}$ is a closed
(compact and without boundary) Riemannian manifold.
In this framework, the system models the coexistence of three immiscible
phases, and the potential $W\in C(\mathbb{R}^2,\mathbb{R}_{\ge 0})$ exhibits a triple-well structure, which means that it admits three
nondegenerate global minima at the same level and satisfies conditions \eqref{eq:W0}, \eqref{eq:W1}, \eqref{eq:W2}, and \eqref{eq:W3}.
We denote the set of minima by $\mathcal{Z}=\{\mathbf{z}_0 = 0,\mathbf{z}_1,\mathbf{z}_2\}\subset\mathbb{R}^2$.
Once again, a $\Gamma$-convergence result applies
(cf.~\cite{MR1051228}), showing that the Allen--Cahn functional provides
a relaxation of the multi-perimeter functional associated with clusters.
More precisely, if $(u_\varepsilon)_{\varepsilon>0}
\subset W^{1,2}({M},\mathbb{R}^2)$ has uniformly bounded Allen--Cahn energy,
then, up to a subsequence, it $L^1$-converges to a function
$u_0 \in L^1({M},\mathbb{R}^2)$ such that
\[
u_0(x) \in \{\mathbf{z}_0,\mathbf{z}_1,\mathbf{z}_2\}
\quad \text{for a.e. } x\in{M}.
\]
The regions $u_0^{-1}(\mathbf{z}_1)$ and $u_0^{-1}(\mathbf{z}_2)$ then form a two-cluster
$\Omega=\{\Omega^1,\Omega^2\}$, while by convention $u_0^{-1}(\mathbf{z}_0)$ is regarded as the exterior chamber.
Moreover, one has the lower bound
\[
\liminf_{\varepsilon\to0}\mathcal{AC}_\varepsilon(u_\varepsilon)
\ge \mathcal{P}_g^\omega(\Omega)
:= \sum_{0\le i<j\le2}
\omega_{ij}\,
\mathcal{H}_g^{N-1}\big(\partial^*\Omega_i \cap \partial^*\Omega_j\big),
\]
where the coefficients $\omega_{ij}>0$ depend only on the triple-well potential $W$.
In addition, they are assumed to satisfy the immiscibility condition \eqref{eq:W0}.

For a small volume regime of the two considered phases (regarding the third one corresponding to $\mathbf{z}_0 = 0$ as the exterior phase), we recall that we study the Allen--Cahn functional under the following constraint
\[
    W^{1,2}_\mathrm{v}({M},\mathbb{R}^2) = \Big\{u \in W^{1,2}({M},\mathbb{R}^2):     \int_{M} u\,  \mathrm{d}v_g = \mathrm{v} \in      \mathbb{R}^2_{> 0}\Big\}.
\]
In this setting, the main difficulty lies in the construction of a suitable
photography map.
Since any Riemannian metric is locally close to the Euclidean one, a natural candidate is obtained by projecting, via the exponential map, a model
two-cluster of volume $\mathrm{v}\in\mathbb{R}^2_{>0}$ defined in a tangent space
$T_p{M} \simeq \mathbb{R}^N$.
In Euclidean space, such a two-cluster is given by the double bubble, which is known to be the unique (up to translations and rotations) isoperimetric cluster and admits an explicit description.

However, a continuous global choice of these projected double bubbles is not
always possible, since their geometry naturally selects a preferred direction
for the relative position of the two chambers.
For this reason, in~\cite{MR4701348} the manifold ${M}$ is assumed to be
\emph{parallelizable}, allowing for a continuous choice of an orthonormal frame 
available on the whole manifold.
This allows one to define, for each $p\in{M}$, a two-cluster contained in a
small ball centered at this point, in a continuous way.
A further issue arises from the fact that, under the exponential map, the
projection of a Euclidean double bubble generally produces a cluster whose
volumes differ from the prescribed vector $\mathrm{v}\in\mathbb{R}^2$, due to curvature distortions.
This difficulty is overcome by exploiting the explicit characterization of the
Euclidean double bubble.
Indeed, since the map $\mathrm{v} \mapsto \Omega_\mathrm{v}$ is well-defined and continuous (up to
translations and rotations), one can compensate for the geometric distortion by
selecting, at each point $p\in M$, a slightly modified two-cluster $\Omega_{\mathrm{v}'(p)}\in\mathcal{C}^1_g(M)$ in the
tangent space, with $\mathrm{v}'(p)\approx \mathrm{v}$ chosen in a continuous manner.

In this way, the resulting clusters on the manifold have exactly the prescribed
volume, and a recovery sequence construction yields the desired photography
map.
Moreover, since for small volumes the clusters are supported in small balls
where the metric tensor is uniformly bi-Lipschitz and close to the identity,
the corresponding energy level $\mathcal{AC}_\varepsilon^{c(\mathrm{v})}$ can be estimated by means of the perimeter
formula for the Euclidean double bubble,
and its first order expansion is again the isoperimetric profile (for vectorial volume), thus leading to an estimate as in~\eqref{eq:cvtoIv}.
Thanks to this estimate, one obtains a \emph{concentration property} for the functions in the sublevel $\mathcal{AC}_\varepsilon^{c(\mathrm{v})}$,
thus the barycenter map is defined, again by using 
an extrinsic barycenter via Nash-embedding.
This whole procedure then gives the following result.
\begin{theorem}[{cf.~\cite[Theorem 1]{MR4701348}}]
\label{thm:vectorial-summary}
    Let $N,m\in\mathbb{Z}_{>0}$ with $N\ge 2$, and $({M},g)$ be a closed parallelizable $N$-dimensional Riemannian manifold.
    There exists $\mathrm{v}^*>0$ such that for every
    $\mathrm{v}=(\mathrm{v}_1,\mathrm{v}_2)\in(0,\mathrm{v}^*)^2$,
    there exists $\varepsilon^*=\varepsilon^*(M,g,\mathrm{v})>0$ satisfying that
    for every $\varepsilon\in(0,\varepsilon^*)$, Eq. \eqref{eq:AC-PDE} with $m=2$ admits at least $\mathrm{cat}({M})+1$ distinct solutions.
    Moreover, in the nondegenerate setting, one can find at least $2\mathcal{P}_1({M})+1$ distinct solutions.
\end{theorem}

%====================================================
\subsection{Multi-bubble case}
%====================================================
For three or more phases, the situation changes dramatically.
The absence of a complete classification of multi-bubble minimizers
leads to serious obstructions, even in the plane.
While partial results are available in special geometries or under
restrictive assumptions, a general theory paralleling the scalar case is currently out of reach.
Below, we discuss some aspects of the proof of Theorem~\ref{thm:intro-vectorial_clusters}.

In~\cite{andrade2024clusters}, a strategy to overcome these difficulties
was introduced, allowing for obtaining multiplicity results for an arbitrary number
of phases without relying on a classification of multi-bubble minimizers.
Nevertheless, two fundamental properties are still available.
As proved by Almgren in his classical work~\cite[Chapter~VI]{MR420406}
(see also~\cite{MR2976521} for a modern exposition),
minimizers of the multi-isoperimetric problem always exist
and are contained in a bounded region
(in particular, no chamber can develop arbitrarily long and thin tentacles).
Moreover, \emph{small} variations of the prescribed volumes can be performed at the cost of only a controlled increase of the perimeter functional (namely, linear with respect to the change of volume).
This procedure is called \emph{volume-fixing lemma}
(cf.~\cite[Theorem 29.14]{MR2976521}).

These properties, however, are not sufficient to produce a continuous selection of isoperimetric clusters with respect to the volume constraint.
What they provide is a weaker form of stability:
loosely speaking, one can construct a \emph{locally} continuous selection of
\emph{almost} isoperimetric clusters in a precise quantitative sense.
This proves to be enough to implement the previous constructions
and obtain a multiplicity result via the photography method.
Indeed, as in the earlier applications, the configurations produced by the
photography map are already almost minimizers.
As a consequence, there is no need to project
(via the exponential map) exact isoperimetric minimizers;
working with almost minimizers is sufficient,
provided that one has suitable control of the perimeter deficit.
Therefore, if the admissible volume variations for which a locally continuous
selection of almost minimizing clusters is available are large enough
to compensate for the volume distortion induced by the curvature of the manifold,
the above machinery applies and a multiplicity result follows.

However, in order to make this strategy work, a structural restriction was adopted.
Indeed, one cannot impose a bound on the \emph{total} volume of the cluster
that guarantees the validity of the construction for all admissible configurations
below such a threshold.
The reason is that, even among clusters with fixed total volume,
one (or more) chambers may have an arbitrarily small volume.
This obstruction comes from the fact that the volume-fixing lemma
depends on the volumes of the individual chambers.
If one chamber becomes too small, the corresponding estimates degenerate,
independently of the total volume of the cluster.
Hence, in the multi-bubble regime, a bound on the total volume alone
does not provide the required control.

For this reason, a different constraint must be imposed,
ensuring that the volume ratios between the chambers remain fixed.
This can be achieved by exploiting the scaling invariance
of the problem in the Euclidean setting.
More precisely, let $\Omega_{\mathrm{v}}$ be an $m$-cluster in $\mathbb{R}^N$
with volume vector $\mathrm{v} \in \mathbb{R}^m_{>0}$.
For $\alpha>0$, we denote by $\alpha \Omega_{\mathrm{v}}$
its rescaled version with volume $\alpha \mathrm{v}$, namely
\[
    x \in (\alpha\Omega)^i 
    \iff
    \frac{x}{\alpha^{1/N}} \in \Omega^i,
    \qquad i = 1,\dots,m.
\]
With this definition, the Euclidean multi-perimeter satisfies
\[
\mathcal{P}^\omega_\delta (\alpha\Omega)
= \alpha^{\frac{N-1}{N}}\,\mathcal{P}^\omega_\delta (\Omega),
\qquad \forall \alpha > 0.
\]

Assume now that the volume-fixing procedure for the minimizer
$\Omega_{\mathrm{v}}$ allows variations of the volume vector
within a ball $B(\mathrm{v},\zeta) \subset \mathbb{R}^m_{>0}$,
for some $\zeta>0$.
In other words, for any $\tilde{\mathrm{v}} \in B(\mathrm{v},\zeta)$
there exists an $m$-cluster $\Omega_{\tilde{\mathrm{v}}}$
with the prescribed volume, and this choice depends continuously
on the volume parameter.
By scaling invariance, the same construction applies to
$\Omega_{\alpha\mathrm{v}}$,
yielding admissible volume variations in the ball
$B(\alpha \mathrm{v},\alpha^{1/N}\zeta)$.
Moreover, the corresponding multi-perimeter deficit
is controlled through the rescaling formula.

This observation is the key point of the construction.
While the prescribed volume rescales like $\alpha$,
the volume distortion induced by projecting the cluster
onto the manifold is of a higher-order, namely of order
$\alpha^{1+2/N}$.
As a consequence, for $0<\alpha\ll1$ sufficiently small
(depending on the fixed vector $\mathrm{v}\in\mathbb{R}^m_{>0}$),
the admissible range of volume variations for
$\Omega_{\alpha\mathrm{v}}\subset M$ is large enough to compensate
for the distortion produced by the curvature.
This allows one to construct, at every point $p\in M$,
a cluster with the correct volume $\alpha \mathrm{v}$,
depending continuously on $p$.
We remark that this construction relies crucially on the compactness
of the manifold $M$ and on the local bi-Lipschitz property
of the exponential map.

In this fashion, we may state the following multiplicity result:
\begin{theorem}[{cf.~\cite[Theorem~A]{andrade2024clusters}}]
	\label{theoremfrombubblestoclusters}
    Let $N,m\in\mathbb{Z}_{>0}$ with $N\ge 2$, and $({M},g)$ be a closed parallelizable $N$-dimensional Riemannian manifold.
    For every $\mathrm{v}\in\mathbb{R}^m_{> 0}$, there exists $\alpha^* = \alpha^*(M,g,\mathrm{v}) > 0$
	such that for every $\alpha \in (0,\alpha^*)$
	there exists $\varepsilon^* = \varepsilon^*(\alpha,\mathrm{v}) > 0$
	such that for every $\varepsilon \in (0,\varepsilon^*)$, Eq. \eqref{eq:AC-PDE} with $m\geq 3$ and $\tilde{\mathrm{v}}=\alpha\mathrm{v}$ admits at least $\mathrm{cat}({M})+1$ distinct solutions.
    Moreover, in the nondegenerate setting, one can find at least $2\mathcal{P}_1({M})+1$ distinct solutions.
\end{theorem}

%====================================================
%====================================================
\section{Open directions}\label{sec:open}
%====================================================
%====================================================
The results surveyed so far show that the photography method provides a robust
variational-topological framework for Allen--Cahn equations and systems in the
singular perturbation regime, both in scalar and vectorial settings, and in the
presence of boundaries. At the same time, several natural extensions remain
largely unexplored or only partially understood. In this section, we collect a
selection of open problems that are closely connected to the results stated in
Theorems~\ref{thm:intro-scalar-closed}--\ref{thm:intro-vectorial_clusters} and that point toward genuinely new analytical and geometric
phenomena.

Rather than aiming for a comprehensive list, we focus on four directions where
the interaction between $\Gamma$-convergence, geometry, and topology raises
new difficulties beyond those already addressed in the scalar and double-bubble
regimes.

%====================================================
\subsection{Clusters on manifolds with boundary}
%====================================================
As we showed before, in the \emph{scalar} theory, Dirichlet and Neumann
boundary conditions lead to markedly different singular limits and multiplicity
patterns, detecting either the topology of $M$ or that of $\partial M$. By
contrast, Theorems~\ref{thm:intro-vectorial_bubble} and~\ref{thm:intro-vectorial_clusters} treat only the \emph{boundaryless} vectorial case,
and no analog is currently available when $\partial M\neq\emptyset$.

From the sharp-interface viewpoint, one expects the $\Gamma$-limit of vectorial
Allen--Cahn systems with boundary to be a multi-isoperimetric problem with
boundary conditions imposed on clusters. However, even the basic geometric
properties of such minimizers, like localization, regularity near $\partial M$, and interaction between junctions and the boundary, are largely unknown,
especially for the case of three or more chambers.
Partial evidence of the geometric complexity of multiphase clusters with boundary can be found in the recent work~\cite{novaga2025}, where isoperimetric triple partitions arising as blow-down limits naturally lead to partition problems in half-spaces, and non-standard configurations are shown to exist in high dimension.
This lack of control obstructs the extension of the photography method,
whose key steps rely on uniform diameter bounds and concentration of small clusters,
and motivates the open problem of extending
Theorems~\ref{thm:intro-vectorial_bubble} and~\ref{thm:intro-vectorial_clusters}
to manifolds with boundary, obtaining multiplicity results for vectorial
Allen--Cahn systems with Dirichlet or Neumann boundary conditions.

%====================================================
\subsection{Nonlocal Allen--Cahn systems}
%====================================================
All results in Theorems~\ref{thm:intro-scalar-closed}--\ref{thm:intro-vectorial_clusters} rely on the local Modica--Mortola $\Gamma$-convergence
to the (multi-)perimeter functional. In nonlocal or fractional Allen--Cahn
models, the Dirichlet energy is replaced by long-range interactions, and the
singular limit is given by a nonlocal perimeter. While $\Gamma$-convergence
results are available in several settings, the geometric and topological
structure of low-energy configurations differs substantially from the local
case.

Nonlocal variants of the Allen--Cahn functional, involving fractional
diffusion operators have been extensively
studied in recent years. In the scalar case, these models admit a sharp--interface
limit governed by nonlocal perimeters and fractional minimal surfaces, together with
existence and qualitative properties of minimizers and their relation with nonlocal CMC surfaces
\cite{MR2675483,MR2948285}.
Multiplicity results for vectorial and genuinely multiphase nonlocal models,
closely related to cluster formation, are far less understood and remain largely
open beyond existence and pattern selection results \cite{MR3614671}.
We also refer to the recent resolution of Yau's conjecture for nonlocal CMC hypersurfaces (cf. \cite{arXiv:2306.07100}), which is still open in the local case.
This suggests that richer behavior is expected for the nonlocal Allen--Cahn systems.

In particular, the concentration and localization mechanisms that underlie the
photography method may fail, or may require substantial
modification, due to long-range effects. An open question is whether analog
of the multiplicity results discussed here persist for
nonlocal Allen--Cahn models, and how the detected topology depends on the
interaction length scale.

%====================================================
\subsection{Diffuse CMCs on noncompact manifolds}
%====================================================
The results presented in this survey are formulated on compact manifolds, where
Palais--Smale compactness and sharp-interface localization can be enforced by
small-volume arguments. On noncompact manifolds, however, minimizing and critical
sequences may drift to infinity, and the direct application of the photography
method breaks down.

From the $\Gamma$-convergence viewpoint, the singular limit is expected to be a
relative isoperimetric problem, possibly localized near a compact core or along
an end of the manifold.
Such sharp-interface problems have been extensively
studied on noncompact spaces under various geometric assumptions
(see~\cite{MR4866981} for a recent survey on this topic).
By contrast, much less is known at the diffuse-interface level: at present,
the Allen--Cahn functional and its minimizers have been analyzed essentially only in the Euclidean setting, where the arguments rely heavily on the presence of a cocompact
action by isometries
(see~\cite{bonforte2024,MR4761254}).

An important open problem is therefore to determine whether versions of
Theorems~\ref{thm:intro-scalar-closed} and~\ref{thm:intro-vectorial_bubble} can be
recovered on noncompact manifolds under suitable geometric assumptions (for
instance, asymptotically flat or cylindrical metrics), possibly by combining the
photography method with concentration-compactness techniques.

%====================================================
\subsection{From interfaces to higher-dimensional vortices}
%====================================================
Theorems~\ref{thm:intro-scalar-closed} and~\ref{thm:intro-vectorial_clusters} crucially rely on the fact that the potential $W$ has a finite
set of nondegenerate minima, so that the singular limit produces codimension-one
interfaces. A different regime arises when the zero set of
the potential $\mathcal{Z}\subset \mathbb{R}^m$ is a smooth manifold.
An important example is provided by Ginzburg--Landau-type models, where the
set of minima of the potential is a smooth manifold (typically $\mathbb S^1 \subset \mathbb{R}^2$)
and the singular limit is characterized by the formation of vortices, namely
topological defects of codimension two, rather than codimension-one interfaces.

A first result in this direction has appeared very recently (see~\cite{corona2025GL}), showing that
photography-type constructions can indeed be implemented in the scalar
Ginzburg--Landau setting to produce topological multiplicity of critical points.
This suggests that the underlying philosophy of the photography method is not
intrinsically tied to codimension-one interfaces, but may extend to regimes
where the singular limit involves higher-codimensional defects.
At the same time, many questions remain open concerning the flexibility of this
approach under different constraints or boundary conditions, and a systematic theory in this direction is still largely to be
developed.

\section*{Declarations}

\noindent{\bf Funding.}
This work was partially supported by Funda\c c\~ao de Amparo \`a Pesquisa do Estado de S\~ao Paulo (FAPESP) and Conselho Nacional de Desenvolvimento Cient\'ifico e Tecnol\'ogico (CNPq). 
J.H.A. acknowledges financial support from FAPESP \#2023/15567-8 and CNPq \#409764/2023-0, \#443594/2023-6, \#441922/2023-6, \#306014/2025-4. 
S.N. acknowledges financial support from FAPESP \#2021/05256-0, \#2023/08246-0 and CNPq \#12327/2021-8, \#441922/2023-6.
R.P. acknowledges financial support from FAPESP \#2023/07697-9, \#2025/07027-9.
\medskip

\noindent{\bf Conflict of interest.}
The authors have no relevant financial or non-financial interests to disclose.
\medskip

\noindent{\bf Data availability.}
Data sharing is not applicable as no datasets were generated or analyzed during the study.
\medskip

\noindent{\bf Ethics approval.}
Not applicable.

\begingroup
\renewcommand*{\bibfont}{\small}
\printbibliography
\endgroup

\end{document}